\documentclass[11pt,a4paper]{article}

\usepackage{textcomp}
\usepackage{a4wide}
\usepackage{graphicx}
\usepackage{authblk}
\usepackage{amssymb}
\usepackage{epstopdf}
\usepackage[sans]{dsfont}
\usepackage[english]{babel}
\usepackage{latexsym}
\usepackage{subfigure} 
\usepackage{color}
\usepackage{float}
\usepackage{amsmath}
\usepackage{amsfonts}
\numberwithin{equation}{section}
\usepackage{amsthm}
\usepackage[numbers]{natbib}
\usepackage[bookmarks=true,colorlinks=true,linkcolor={blue},urlcolor={blue}, citecolor={blue},pdfstartview={XYZ null null 1.22}]{hyperref}

\newtheorem{remark}{Remark}
\newtheorem{proposition}{Proposition}
\newtheorem{theorem}{Theorem}

\newtheorem{lemma}{Lemma}
\newtheorem{algorithm}{{Algorithm}}

\def\RR{\mathbb R}

\def\be{\begin{equation}}
\def\ee{\end{equation}}
\def\bea{\begin{eqnarray}}
\def\eea{\end{eqnarray}}

\begin{document}
\title{Uncertainty quantification for kinetic models \\ in socio--economic and life sciences}

\author{Giacomo Dimarco\thanks{Department of Mathematics and Computer Science, University of Ferrara, Via Machiavelli 35, 44121 Ferrara, Italy ({\tt giacomo.dimarco@unife.it}).} \and Lorenzo Pareschi\thanks{Department of Mathematics and Computer Science, University of Ferrara, Via Machiavelli 35, 44121 Ferrara, Italy ({\tt lorenzo.pareschi@unife.it}).} \and
Mattia Zanella\thanks{Department of Mathematical Sciences, Politecnico di Torino, Corso Duca degli Abruzzi 24, Torino, Italy ({\tt mattia.zanella@polito.it}).}
}
\maketitle

\begin{abstract}
Kinetic equations play a major rule in modeling large systems of interacting particles. Recently the legacy of classical kinetic theory found novel applications in socio-economic and life sciences, where processes characterized by large groups of agents exhibit spontaneous emergence of social structures. Well-known examples are the formation of clusters in opinion dynamics, the appearance of inequalities in wealth distributions, flocking and milling behaviors in swarming models, synchronization phenomena in biological systems and lane formation in pedestrian traffic. The construction of kinetic models describing the above processes, however, has to face the difficulty of the lack of fundamental principles since physical forces are replaced by empirical social forces. These empirical forces are typically constructed with the aim to reproduce qualitatively the observed system behaviors, like the emergence of social structures, and are at best known in terms of statistical information of the modeling parameters. For this reason the presence of random inputs characterizing the parameters uncertainty should be considered as an essential feature in the modeling process. In this survey we introduce several examples of such kinetic models, that are mathematically described by nonlinear Vlasov and Fokker--Planck equations, and present different numerical approaches for uncertainty quantification which preserve the main features of the kinetic solution.  
\end{abstract}

\tableofcontents

\section{Introduction}
Kinetic models describing the collective behavior of a large group of interacting agents have attracted a lot of interest in the recent years in view of their potential applications to various fields, like sociology, economy, finance and biology \cite{AH,AjBE,AHP,APTZ,APZc,CChoH,CFRT,CFTV,CPT,CristPiccTos01,CristPiccTos02,DLR,DueWolf,NPT,PT2}. One of the major difficulties in applying the classical toolbox of kinetic theory to these systems is the lack of fundamental principles which define the microscopic dynamic. In addition, experimental results are typically non reproducible and, as a consequence, the model construction is dictated by its ability to describe qualitatively the system behavior and the formation of emergent social structures. A degree of uncertainty is therefore implicitly embedded in such models, since most modeling parameters can be assigned only as statistical information from experimental results \cite{albi2015MPE,Ball_etal,BFHM,DMPW,KTIHC}.         

From a mathematical viewpoint, the kinetic models we will consider in the present survey are characterized by nonlinear Vlasov--Fokker--Planck equations with random inputs taking into account uncertainties in the initial data, in the interaction terms and/or in the boundary conditions. The models describe the evolution of a  distribution function $f=f(\theta,x,w,t)$, $t\ge 0$, $x\in\RR^{d_x}$, $w\in\RR^{d_w}$, $d_x,d_w\ge 1$, and $\theta\in\Omega\subseteq\RR^{d_\theta}$ \emph{a random field}, accordingly to
\be\begin{split}\label{eq:MF_general}
&\partial_t f+\mathcal L[f] = \nabla_w \cdot \left[ \mathcal B[f]f+\nabla_w (Df) \right],
\end{split}\ee
where  $\mathcal L[\cdot]$ is a linear operator describing the agents' dynamics with respect to the $x-$variable,{typically $\mathcal L[f]=w\cdot\nabla_x f$,} $\mathcal B[\cdot]$ is a non--local operator of the form 
\be
\label{eq:Bf}
\mathcal B[f](\theta,x,w,t) = \int_{\RR^{d_x}}\int_{\RR^{d_w}}P(x,x_*;w,w_*,\theta)(w-w_*)f(\theta,x_*,w_*,t)dw_*dx_*,
\ee
and $D(\theta,w)\ge 0$, for all $w\in\RR^{d_w}$, is a function describing the local relevance of the diffusion. {We refer to \cite{Cer,DPR,Risk,Vill,Vlas} for an introduction to the subject in relation with kinetic theory.} In the rest of the chapter, to avoid unnecessary difficulties, we will mainly restrict to the case of a one-dimensional random input $d_{\theta}=1$ distributed as $p(\theta)$.
In the homogeneous case $f=f(\theta,w,t)$, $\mathcal L[f]\equiv 0$ the kinetic models are characterized by nonlinear Fokker-Planck equations. 

\paragraph{\it The classic Fokker-Planck equation with uncertainties}
The most classical example is represented by the linear Fokker-Planck model obtained for $P\equiv 1$ corresponding to 
\be
\mathcal B[f](\theta,w) = (w-u(\theta)),\qquad  D(\theta)=T(\theta),
\ee  
where
\[
u(\theta)=\int_{\RR^{d_w}} f(\theta,w,t)w\,dw,\quad T(\theta)=\frac1{d_{w}}\int_{\RR^{d_w}} f(\theta,w,t) (w-u(\theta))^2 f(\theta,w,t)\,dw
\]
are the (conserved) mean velocity and the temperature of the particles. In the above expressions we assumed an uncertain initial data such that $\int_{\RR^{d_w}} f(\theta,w,0)\,dw=1$ for all $\theta\in\Omega$. The stationary solution in this case is represented by a Maxwellian distribution with uncertain momentum and temperature given by
\be 
\label{eq:maxw}
f^{\infty}(\theta,w)=\frac{1}{(2\pi T(\theta))^{d_w/2}}\exp\left\{-\frac{|w-u(\theta)|^2}{2T(\theta)}\right\}.
\ee

\paragraph{\it Opinion formation with uncertain interaction}
A kinetic Fokker-Planck model of opinion formation for $w\in I=[-1,1]$, where $\pm 1$  denote the two extremal opinions, corresponds to the choices \cite{PT2,T}
\be\label{eq:BD_opinion}
\mathcal B[f](\theta,w,t) = \int_{I} P(\theta,w,w_*)(w-w_*)f(\theta,w_*,t)dw_*,\quad D(w) =\dfrac{\sigma^2}{2} (1-w^2)^2.
\ee
In the above nonlocal interaction term $P(\theta,\cdot,\cdot)\in [0,1]$ is a function taking into account uncertainties in the compromise propensity between the agents' opinions. 

In the simple case $P(w,w_*,\theta)=P(\theta)$ and deterministic initial data, the model preserves the mean opinion $u=\int_{I}wf(\theta,w,t)dw$ and we can analytically compute the steady state distribution 
\be\label{eq:exact_ss_op}
f^{\infty}(\theta,w) = \dfrac{C}{(1-w^2)^2}\left( 1+w \right)^{\frac{P(\theta)u}{2\sigma^2}}(1-w)^{\frac{P(\theta)u}{2\sigma^2}} \exp\Big\{ -\dfrac{P(\theta)(1-uw)}{\sigma^2(1-w^2)} \Big\},
\ee
with $C>0$ a normalization constant. 

\paragraph{\it Wealth distribution with uncertain diffusion}
If we now consider $w\in [0,\infty)$ a measure of the agents' wealth, a Fokker-Planck model describing the wealth evolution of agents is obtained taking \cite{CPT,PT2} 
\be\label{eq:BD_wealth}
\mathcal B[f](\theta,w,t) = \int_{[0,\infty]} a(w,w_*)(w-w_*)f(\theta,w_*,t)dw_*,\quad D(\theta,w) = \frac{\sigma(\theta)^2}{2}w^2,
\ee
where the term $\sigma(\theta)$ characterizes the uncertain strength of diffusion. 
An explicit expression of the steady state distribution is given in the case $a(w,w_*)\equiv 1$ 
\be\label{eq:exact_ss_we}
f^{\infty}(\theta,w) = \dfrac{(\mu(\theta)-1)^{\mu(\theta)}}{\Gamma(\mu(\theta))w^{1+\mu(\theta)}}\exp\Big\{ -\dfrac{\mu(\theta)-1}{w} \Big\},
\ee
where $\Gamma(\cdot)$ is the Gamma function and $\mu(\theta)=1+2/\sigma^2(\theta)$ is the so--called Pareto exponent, which is now dependent on the random input.

\paragraph{\it Swarming models with uncertainties}
As a final example we consider a kinetic model for the swarming behavior \cite{BS2012,CChoH,CFRT,CFTV,CHL,CS,DegLMP,DOCBC,HLL,HT}. In particular we focus on a model with self--propulsion and uncertain diffusion, see \cite{BarDeg,BCCD}. The dynamics for the density $f=f(\theta,x,w,t)$ of agents in position $x\in\RR^{d_x}$ with velocity $w\in\RR^{d_w}$ is described by the Vlasov-Fokker-Planck equation (\ref{eq:MF_general}) characterized by 
\be
\begin{split}
\label{eq:swarming_general}
\mathcal L[f]=w\cdot\nabla_x f, \qquad \mathcal B[f](\theta,x,w,t) = \alpha w(1-|w|^2)+(w-u_f(\theta,x,t)),
\end{split}
\ee
where 
\be
u_f(\theta,x,t) = \dfrac{\int_{\RR^{d_x}\times\RR^{d_w}}K(x,y)wf(\theta,y,w,t)\,dw\,dy}{\int_{\RR^{d_x}\times \RR^{d_w}}K(x,y)f(\theta,y,w,t)\,dw\,dy},
\ee
with $K(x,y)>0$ a localization kernel, $\alpha>0$ a self--propulsion term and $D(\theta)>0$ the uncertain noise intensity. 

In the space--homogeneous case $f=f(\theta,w,t)$, stationary solutions have the form 
\be
f^{\infty}(w,\theta) = C\exp \left\{-\dfrac{1}{D(\theta)}\left(\alpha\dfrac{|w|^4}{4}+(1-\alpha)\dfrac{|w|^2}{2}-u_{f^{\infty}}(\theta)\cdot w \right) \right\},
\ee
with $C>0$ a normalization constant and $$u_{f^{\infty}}(\theta) = \dfrac{\int_{\RR^{d_w}}wf^{\infty}(w,\theta)dw}{\int_{\RR^{d_w}}f^{\infty}(w,\theta)dw}.$$
{We stress that, in all the above reported examples, uncertainty may be present in other modeling parameters by further increasing the dimensionality and the complexity of the kinetic model. }

The development of numerical methods for kinetic equations presents several difficulties due to the high dimensionality and the intrinsic structural properties of the solution. Non negativity of the distribution function, conservation of invariant quantities, entropy dissipation and steady states are essential in order to compute qualitatively correct solutions. Preservation of these structural properties is even more challenging in presence of uncertainties which contribute to increase the dimensionality of the problem. We refer to \cite{DP15,JinParma,Son} for recent surveys on numerical methods for kinetic equations in the deterministic case. 

For this reason we will focus on the construction of numerical methods for uncertainty quantification (UQ) which preserves the structural properties of the kinetic equation and, in particular, which are able to capture the correct steady state of the problem with arbitrary accuracy. We will discuss different numerical approaches based on the major techniques used for uncertainty quantification. In the deterministic case, similar approaches for nonlinear Fokker-Planck equations were previously derived in \cite{BuetDellacherie2010,BCDS,ChCo,LLPS, MB,SG}. Related methods for the case of nonlinear degenerate diffusion equations were proposed in \cite{BCF,CJS} and with nonlocal terms in \cite{BCW,CCH}. We refer also to \cite{AlbiPareschi2013ab} for the development of methods based on stochastic approximations and to \cite{Gosse} for a recent survey on schemes which preserve steady states of balance laws and related problems.

The simplest class of methods for quantifying uncertainty in partial differential equations (PDEs) are the stochastic collocation methods. Stochastic collocation methods are non-intrusive, so they preserve all properties of the deterministic numerical scheme, and easy to parallelize. In Section 3 we describe the structure preserving methods recently developed in \cite{PZ1, PZ2} together with a collocation approach and show how the resulting schemes preserve non negativity, conservation and entropy dissipation. In addition they capture the steady states with arbitrary accuracy and may achieve high convergence rates (spectral convergence for smooth solutions). Next in Section 4, we consider the closely related class of statistical sampling methods, most notably Monte Carlo (MC) sampling. In order to address the slow convergence of MC methods, we discuss here the development of Monte Carlo methods based on a Micro--Macro decomposition approach {introduced in \cite{DP}.} These methods preserve the structural properties of the kinetic problem, are capable to significantly reduce the statistical fluctuations of standard Monte Carlo and increase their computational efficiency by reducing the number of statistical samples in time. 
Section 5 is devoted to stochastic Galerkin methods based on generalized Polynomial Chaos (gPC). Although these deterministic methods may achieve high convergence rates for smooth solutions, they suffer from the disadvantage that they are highly intrusive and that increase the computational complexity of the problem. As a consequence, the main physical properties of the solution are typically lost at a numerical level. For this class of methods we show how to construct generalized Polynomial Chaos schemes based on the Micro--Macro formalism which preserve the steady states of the system \cite{PZ3}. Finally, in Section 6 several numerical applications to problem in socio-economy and life sciences are presented. 

\section{Preliminaries}
{In this section we recall some analytical properties of the considered kinetic models which will be useful for the development of the different numerical methods.} 
Except in some simple case, a precise analytic description of the global equilibria of equation  \eqref{eq:MF_general} is very difficult \cite{CMV,CarTos,ToVil}. A deeper insight into the large time behavior can be achieved by resorting to the asymptotic behavior of the corresponding space homogeneous models, leading to {nonlinear} Fokker--Planck type equations \cite{Tos1999}.

\subsection{Fokker-Planck type equations}
In {the space homogeneous} case the distribution function reduces to $f=f(\theta,w,t)$, $w\in\RR^{d_w}$, $\theta\in\RR^{d_\theta}$, $t>0$ and is solution of the following problem
\be\label{eq:FP_general}
\partial_t f(\theta,w,t) = \mathcal J(f,f)(\theta,w,t) ,
\ee
where
\be\label{eq:J_def}
\mathcal J(f,f)(\theta,w,t) = \nabla_w\cdot \Big[ \mathcal B[f](\theta,w,t)f(\theta,w,t)+\nabla_w D(\theta,w)f(\theta,w,t) \Big], 
\ee
together with an initial datum $f(\theta,w,0) = f_0(\theta,w)$ and suitable boundary conditions on $w\in\RR^{d_w}$. 

We review in the present stochastic setting the classical results for the trend to equilibrium of the problem \eqref{eq:FP_general} in the simplified case of a one-dimensional problem $w\in I\subseteq \RR$ with a linear drift term, i.e.
\be\label{eq:FP_linearcase}
\partial_t f(\theta,w,t) = \partial_w \Big[ (w-u)f(\theta,w,t)+\partial_w ({D(\theta,w)}f(\theta,w,t))\Big].
\ee
Conservation of mass is imposed on the previous equation by considering suitable boundary conditions \cite{PT2}. The stochastic stationary solution $f^{\infty}(\theta,w)$ of equation \eqref{eq:FP_linearcase} is given by the solution of 
\[
(w-u)f^{\infty}(\theta,w)+\partial_w D(\theta,w)f^{\infty}(\theta,w).
\]
The stochastic Fokker--Planck equation \eqref{eq:FP_linearcase} may be rewritten in the equivalent forms
\be\label{eq:Landau_log}
\partial_t f(\theta,w,t) = \partial_w \Big[ D(\theta,w)f(\theta,w,t)\partial_w \log \dfrac{f(\theta,w,t)}{f^{\infty}(\theta,w,t)}\Big],
\ee
which corresponds to the \emph{stochastic Landau form}, whereas the stochastic non logarithmic Laundau form of the equation is the following
\be\label{eq:Landau_nonlog}
\partial_t f(\theta,w,t) = \partial_w \Big[ D(\theta,w)f^{\infty}(\theta,w,t)\partial_w \dfrac{f(\theta,w,t)}{f^{\infty}(\theta,w)} \Big].
\ee
Convergence to equilibrium is usually determined through measures of the entropy production. We define the relative entropy for all positive functions $f,\tilde f$ as follows
\be\label{eq:rel_entropy}
\mathcal H[f,\tilde f](\theta,w,t) = \int_I f(\theta,w,t) \log\left(\dfrac{f(\theta,w,t )}{\tilde f(\theta,w,t)} \right)dw,
\ee
we have \cite{FPTT_K}
\be\label{eq:rel_ent_def}
\dfrac{d}{dt}\mathcal H[f,f^{\infty}](\theta,w,t) = -\mathcal I_D[f,f^{\infty}](\theta,w,t),
\ee
where the dissipation functional $\mathcal I_D[\cdot,\cdot]$ is defined as 
\be\begin{split}
\mathcal I_D[f,f^{\infty}]& = 
 \int_{\mathcal I} D(\theta,w)f(\theta,w,t)\left(\partial_w \log \left( \dfrac{f(\theta,w,t)}{f^{\infty}(\theta,w)}\right)\right)^2dw.
\end{split}\ee
In the classical setting $w\in\RR$ and $D(\theta,w)={T(\theta)}$, where $T(\theta)$ is the temperature, the steady state is given by the Maxwellian density (\ref{eq:maxw}) with $d_w=1$  and relation \eqref{eq:rel_ent_def} coupled with the log--Sobolev inequality  
\be
 \mathcal H[f,f^{\infty}](\theta,w,t)\le \dfrac{1}{2}\mathcal I_D[f,f^{\infty}](\theta,w,t),
\ee
leads to the exponential decay of the relative entropy as proved in the following result \cite{Tos1999}.
\begin{theorem}\label{teo:toscani}
Let $f(\theta,w,t)$ be the solution to the initial value problem
\[
\partial_t f(\theta,w,t) = \partial_w (w-u(\theta))f(\theta,w,t)+T(\theta)\partial_w^2 f(\theta,w,t)
\]
with the initial condition $f(\theta,w,0)=f_0(\theta,w)$ with finite entropy. Then $f(\theta,w,t)$ converges for all $\theta\in\Omega$ to $f^{\infty}(\theta,w)$ given by (\ref{eq:maxw}) and
\[
\mathcal H[f,f^{\infty}]\le e^{-2t/T(\theta)}\mathcal H[f_0,f^{\infty}].
\]
\end{theorem}
For more general diffusion functions ${D(\theta,w)}$ analogous log--Sobolev inequality are not available. A strategy to study the convergence to equilibrium is to investigate the relation of relative entropy with the relative weighted Fisher information, see \cite{CarTos,FPTT_K,MJT,Tos1999} for more details. 

\subsection{Micro--Macro formulation}
\label{sec:micro-macro}
In this paragraph we describe the Micro--Macro approach to kinetic equations of the form \eqref{eq:FP_general}. The approach is based on the classical Micro--Macro decomposition originally developed by Liu and Yu in \cite{LY} for the fluid limit of the Boltzmann equation. This method has been fruitfully employed for the development of numerical methods by several authors (see \cite{BLM,CrLe,DL,CDL,LemMie,Yan} and the references therein). These techniques has been also recently developed in \cite{DQP,FPR,PR} to construct spectral methods for the collisional operator of the Boltzmann equation that preserves exactly the Maxwellian steady state of the system. 
{Since under} suitable regularity assumptions on the initial distribution the Fokker-Planck equation admits a unique steady state solution $f^{\infty}(\theta,w)$, {the Micro--Macro formulation is obtained decomposing the solution of the differential} problem into the equilibrium part $f^{\infty}$ and the non--equilibrium part $g$ as follows
\be\label{eq:f_MM}
f(\theta,w,t) = f^{\infty}(\theta,w)+ g(\theta,w,t),
\ee
where $g(\theta,w,t)$ is a distribution function such that 
\[
\int_{\RR^{d_w}}\phi(w)g(\theta,w,t)dw = 0
\]
for some moments $\phi(w)=1,w$. The above decomposition \eqref{eq:f_MM} applied to the Fokker-Planck problem \eqref{eq:FP_general}--\eqref{eq:J_def} yields the following result.
\begin{proposition}
If the nonlinear Fokker-Planck equations (\ref{eq:FP_general})-(\ref{eq:J_def}) admits the unique equilibrium state $f^{\infty}(\theta,w)$,
the differential operator $\mathcal J(\cdot,\cdot)$ defined in \eqref{eq:J_def} with $\mathcal B[f]$ given by (\ref{eq:Bf}) may be rewritten as
\be
\mathcal J(f,f)(\theta,w,t) = \mathcal J(g,g)(\theta,w,t) + \mathcal N(f^{\infty},g)(\theta,w,t),
\ee
where $\mathcal N(\cdot,\cdot)$ is a linear operator defined as 
\[
\mathcal N(f^{\infty},g)(\theta,w,t) = \nabla_w \Big[\mathcal B[f^{\infty}]g(\theta,w,t)+\mathcal B[g]f^{\infty}(\theta,w)\Big].
\]
 The only admissible steady state solution of the problem 
\begin{equation}
\label{eq:micro-macro}
\begin{cases}
\partial_t g(\theta,w,t) = \mathcal J(g,g)(\theta,w,t)+ \mathcal N(f^{\infty},g)(\theta,w,t),\\
f(\theta,w,t) = f^{\infty}(\theta,w)+g(\theta,w,t)
\end{cases}\end{equation}
is given by $g^{\infty}(\theta,w)\equiv 0$. 
\end{proposition}
The proof is an immediate consequence of the fact that at the steady state we have $\mathcal J(f^{\infty},f^{\infty})=0$. Note that the steady state solution of the reformulated problem (\ref{eq:micro-macro}) is therefore independent of the uncertainty. 

\begin{remark}
Under suitable assumptions, see Theorem \ref{teo:toscani}, {one can show that} $f(\theta,w,t)$ exponentially decays to the equilibrium solution. As a consequence, the non--equilibrium part of the Micro--Macro approximation $g(\theta,w,t)$ exponentially decays to $g^{\infty}(\theta,w)\equiv 0$ for all $\theta\in\Omega$. 
\end{remark}

\section{Collocation methods}
{One of the most popular computational approaches for UQ relies on  
the class of collocation methods \cite{X, XH}. These methods are non intrusive and permit to couple existing solvers for the PDEs without random inputs with techniques for the quantification of the uncertainty. Moreover, since the structure of the solution remains unchanged, numerical analysis of  collocation methods is a straightforward consequence of the results obtained for the underlying method
used for solving the original equation. }

In the following, since the linear transport part in (\ref{eq:MF_general}) can be discretized using standard approaches, see {for instance} \cite{DP15}, we concentrate on homogeneous Fokker-Planck problem of the form \eqref{eq:FP_general}-\eqref{eq:J_def}. 

Collocation methods consist in solving the problem in a finite set of nodes $(\theta_k)_{k=0}^M$ of the random field. In this class of methods belongs the usual Monte Carlo sampling (MC) which will be treated in Section \ref{sec:4}. If the distribution of the random input $\theta\sim p(\theta)$ is known, an efficient way to treat the uncertainty is to select the nodes in the random space according to Gaussian quadrature rules related to such distribution. This is straightforward in the univariate case, whereas becomes more challenging in the multivariate case \cite{Du}. 

For each $k=0,\dots,M$ we obtain a totally deterministic and decoupled problem since the value of the random variable is fixed. Therefore, solving this system of equations poses no difficulty provided one has a well--established deterministic algorithm. The result is an ensemble of $M+1$ deterministic solutions which can be post--processed to recover the statistical values of interest. For example, in the univariate case if $(\omega_k)_{k=0}^M$ are the Gaussian weights on $\Omega\subseteq \RR$ corresponding to $p(\theta)$ we can use the approximations
\begin{eqnarray}
\mathbb{E}[f](w,t)&=&\int_{\Omega} f(\theta,w,t) p(\theta)\,d\theta \approx \mathbb{E}_M [f](w,t)= \sum_{k=0}^M \omega_k f(\theta_k,w,t),\\
\nonumber
Var[f](w,t)&=&\int_{\Omega} (f(\theta,w,t)-\mathbb{E}[f](w,t))^2  p(\theta)\,d\theta \\[-.25cm]
\\[-.25cm]
\nonumber
& \approx&  Var_M[f](w,t)=\sum_{k=0}^M \omega_k (f(\theta_k,w,t)-\mathbb{E}_M [f](w,t))^2,
\end{eqnarray}
where $\mathbb{E}[\cdot]$ and $Var[\cdot]$ denote the mean and the variance respectively.
In the following we concentrates on the construction of numerical schemes which preserve the structural properties of the solution, like non--negativity, entropy dissipation and accurate asymptotic behavior \cite{PZ1, PZ2}. These properties are essential for a correct description of the underlying physical problem. 

\subsection{Structure preserving methods}\label{sec:structure_preserving}
In the one-dimensional case $d_w=1$ for all $k=0,\dots,M$ the Fokker-Planck equation (\ref{eq:FP_general})-\eqref{eq:J_def} may be written as
\be\label{eq:FP_flux}
\partial_t f(\theta_k,w,t) = \partial_w \mathcal F[f](\theta_k,w,t),\quad w\in {I} \subseteq \RR
\ee
where now
\be
\label{eq:flux}
\mathcal F[f](\theta_k,w,t) = ( \mathcal B[f](\theta_k,w,t)+D'(w)) f(\theta_k,w,t)+D(w)\partial_w f(\theta_k,w,t)
\ee
using the compact notation $D'(w)=\partial_w D(w)$. Typically, when ${I}$ is a finite size set the problem is complemented with no-flux boundary conditions at the extremal points. In the sequel we assume $D(w)> 0$ in the internal points of $I$.

We introduce an uniform spatial grid $w_i \in I$, such that $w_{i+1}-w_i=\Delta w$. We denote as usual $w_{i\pm 1/2}=w_i\pm \Delta/2$ and consider a conservative discretization of \eqref{eq:FP_flux} 
\be\label{eq:dflux}
\frac{d}{dt}f_i(\theta_k,t) = \dfrac{\mathcal F_{i+1/2}[f](\theta_k,t)-\mathcal F_{i-1/2}[f](\theta_k,t)}{\Delta w},
\ee
where for each $t\ge 0$ $\mathcal F_{i\pm 1/2}[f](\theta_k,t)$ is the numerical flux function characterizing the discretization. 

Let us set ${\mathcal C}[f](w,\theta_k,t)=\mathcal B[f](w,\theta_k,t)+D'(w)$ and adopt the notations $D_{i+1/2}=D(w_{i+1/2})$, $D'_{i+1/2}=D'(w_{i+1/2})$. We will consider a general flux function which is combination of the grid points $i+1$ and $i$ 
\be\begin{split}\label{eq:CC_flux}
\mathcal F_{i+ 1/2}[f] = \tilde{\mathcal C}^k_{i+1/2}\tilde{f}_{i+1/2}(\theta_k,t)+D_{i+1/2}\dfrac{f_{i+1}(\theta_k,t)-f_i(\theta_k,t)}{\Delta w},
\end{split}\ee
where 
\be\label{eq:f_CC}
\tilde{f}_{i+1/2}(\theta_k,t)=(1-\delta^k_{i+1/2})f_{i+1}(\theta_k,t)+\delta^k_{i+1/2}f_i(\theta_k,t).
\ee
 For example, the standard approach based on central difference is obtained by considering for all $i$ the quantities 
 \[
 \delta_{i+1/2}^k=1/2, \qquad \tilde{\mathcal{C}}^k_{i+1/2}=\tilde{\mathcal C}[f](w_{i+1/2},\theta_k,t).
 \]
 It is well-known, however, that such a discretization method is subject to restrictive conditions over the mesh size $\Delta w$ in order to keep non negativity of the solution.  

Here, we aim at deriving suitable expressions for the family of weight functions $\delta_{i+1/2}^k$ and for $\tilde{\mathcal{C}}^k_{i+1/2}$ in such a way that the method yields nonnegative solutions, without restriction on $\Delta w$, and preserves the steady state of the system with arbitrary order of accuracy.

First, observe that at the steady state the numerical flux should vanish. From \eqref{eq:CC_flux} we get 
\be\label{eq:rapp_1}
\dfrac{f_{i+1}(\theta_k,t)}{f_i(\theta_k,t)} = \dfrac{-\delta_{i+1/2}^k\tilde{\mathcal{C}}^k_{i+1/2}+\dfrac{D_{i+1/2}}{\Delta w}}{(1-\delta_{i+1/2}^k)\tilde{\mathcal{C}}^k_{i+1/2}+\dfrac{D_{i+1/2}}{\Delta w}}.
\ee
Similarly, if we consider the analytical flux imposing ${\cal F}[f](\theta_k,w,t) \equiv 0$, we have
\be\label{eq:FP_steady}
D(w)\partial_w f(\theta_k,w,t) = -(\mathcal B[f](\theta_k,w,t)+D'(w))f(\theta_k,w,t),
\ee
which is in general not solvable, except in some special cases due to the nonlinearity on the right hand side. We may overcome this difficulty in the quasi steady-state approximation integrating equation \eqref{eq:FP_steady} on the cell $[w_i,w_{i+1}]$ 
\be\begin{split}
&\int_{w_i}^{w_{i+1}}\dfrac{1}{f(\theta_k,w,t)}\partial_w f(\theta_k,w,t)dw =\\
&\qquad\qquad\qquad\qquad -\int_{w_i}^{w_{i+1}}\dfrac{1}{D(w)}(\mathcal B[f](\theta_k,w,t)+D'(w))dw,
\end{split}\ee
which gives
\be\label{eq:quasi_SS}
\dfrac{f(\theta_k,w_{i+1},t)}{f(\theta_k,w_i,t)} = \exp \left\{ -\int_{w_i}^{w_{i+1}}\dfrac{1}{D(w)}(\mathcal B[f](\theta_k,w,t)+D'(w))dw  \right\}.
\ee 
Now, by equating the ratio $f_{i+1}(\theta_k,t)/f_i(\theta_k,t)$ and $f(\theta_k,w_{i+1},t)/f(\theta_k,w_i,t)$ in \eqref{eq:rapp_1}--\eqref{eq:quasi_SS} for the numerical and exact flux respectively, and setting
\be\label{eq:B_tilde}
\tilde{\mathcal{C}}^k_{i+1/2}=\dfrac{D_{i+1/2}}{\Delta w}\int_{w_i}^{w_{i+1}}\dfrac{\mathcal B[f](\theta_k,w,t)+D'(w)}{D(w)}dw
\ee
we recover
\be\label{eq:delta}
\delta_{i+1/2}^k = \dfrac{1}{\lambda_{i+1/2}^k}+\dfrac{1}{1-\exp(\lambda_{i+1/2}^k)}, 
\ee
where
\be\label{eq:lambda_high}
\lambda_{i+1/2}^k=\int_{w_i}^{w_{i+1}}\dfrac{\mathcal B[f](\theta_k,w,t)+D'(w)}{D(w)}dw=\frac{\Delta w\,\tilde{\mathcal{C}}_{i+1/2}^k}{D_{i+1/2}}.
\ee
We have the following result \cite{PZ2}
\begin{proposition}
The numerical flux function \eqref{eq:CC_flux}-\eqref{eq:f_CC} with $\tilde{\mathcal{C}}^k_{i+1/2}$ and $\delta_{i+1/2}^k$ defined by \eqref{eq:B_tilde} and \eqref{eq:delta}-\eqref{eq:lambda_high} vanishes when the corresponding flux (\ref{eq:flux}) is equal to zero over the cell $[w_i,w_{i+1}]$. Moreover the nonlinear weight functions $\delta_{i+1/2}^k$ defined by \eqref{eq:delta}-\eqref{eq:lambda_high} are such that $\delta_{i+1/2}^k \in (0,1)$.  
\end{proposition}
By discretizing \eqref{eq:lambda_high} through the midpoint rule 
\be
\int_{w_i}^{w_{i+1}}\dfrac{\mathcal B[f](\theta_k,w,t)+D'(w)}{D(w)}dw = \dfrac{\Delta w(\mathcal B_{i+1/2}(\theta_k,t)+D'_{i+1/2})}{D_{i+1/2}}+O(\Delta w^3),
\ee
we obtain the second order method defined by 
\be\label{eq:lambdamid}
\lambda_{i+1/2}^{k,\textrm{mid}} = \dfrac{\Delta w(\mathcal B_{i+1/2}(\theta_k,t)+ D'_{i+1/2})}{D_{i+1/2}}
\ee
and 
\be
\label{eq:deltamid}
\delta_{i+1/2}^{k,\textrm{mid}} = \dfrac{D_{i+1/2}}{\Delta w(\mathcal B_{i+1/2}(\theta_k,t)+ D'_{i+1/2})}+\dfrac{1}{1-\exp(\lambda_{i+1/2}^{k,\textrm{mid}})}.
\ee
Higher order accuracy of the steady state solution can be obtained using suitable higher order quadrature formulas for the integral \eqref{eq:B_tilde}. We will refer to this type of schemes as structure preserving Chang-Cooper (SP-CC) type schemes.

Some remarks are in order.

\begin{remark}\label{rem:grad}~
\begin{itemize}
\item If we consider the limit case $D_{i+1/2}\to 0$ in \eqref{eq:lambdamid}-\eqref{eq:deltamid} we obtain the weights
\[
\delta^k_{i+1/2}=
\left\{
\begin{array}{cc}
 0, & \mathcal B_{i+1/2}(\theta_k,t)>0,  \\
 1, & \mathcal B_{i+1/2}(\theta_k,t)<0  \\
\end{array}
\right.
\]
and the scheme reduces to a first order upwind scheme for the corresponding aggregation equation.
\item
For linear problems of the form $\mathcal B[f](\theta_k,w,t)=\mathcal B(\theta_k,w)$ the exact stationary state $f^{\infty}(w,\theta_k)$ can be directly computed from the solution of
\be\label{eq:FP_steady2}
D(w)\partial_w f^{\infty}(\theta_k,w) = -(\mathcal B(\theta_k,w)+D'(w))f^{\infty}(\theta_k,w),
\ee
together with the boundary conditions. Explicit examples of stationary states will be reported in the last section. Using the knowledge of the stationary state we have 
\be\begin{split}\label{eq:SS}
\dfrac{f^{\infty}_{i+1}(\theta_k)}{f^{\infty}_i(\theta_k)} &= \exp \left\{ -\int_{w_i}^{w_{i+1}}\dfrac{1}{D(w)}(\mathcal B(\theta_k,w)+D'(w))dw  \right\}\\
&=\exp \left(-\lambda^{\infty}_{i+1/2}(\theta_k) \right),
\end{split}\ee 
therefore 
\be
\label{eq:lambda_inf}
\lambda^{\infty}_{i+1/2}(\theta_k)=\log \left(\frac{f_{i}^{\infty}(\theta_k)}{f_{i+1}^\infty(\theta_k)}\right)\ee
and
\be\label{eq:delta_inf}
\delta^{\infty}_{i+1/2}(\theta_k) = \dfrac{1}{\log(f_i^{\infty}(\theta_k))-\log(f_{i+1}^{\infty}(\theta_k))}+\dfrac{f_{i+1}^{\infty}(\theta_k)}{f_{i+1}^{\infty}(\theta_k)-f_{i}^{\infty}(\theta_k)}. 
\ee
In this case, the numerical scheme preserves the steady state exactly. 
\item The cases of higher dimension $d\ge 2$ may be derived similarly using dimensional splitting (see \cite{PZ1} for details).  
\end{itemize}
\end{remark}

\subsubsection{Main properties}\label{sec:prop_SP}
In the following we recall some results on the preservation of the structural properties, like non negativity and entropy dissipation. 

\paragraph{\em Non negativity}
Concerning non negativity, first we report a result for an explicit time discretization scheme \cite{PZ1}. We introduce a time discretization $t^n=n\Delta t$ with $\Delta t>0$ and $n = 0,\dots,T$ and consider the simple forward Euler method
\be\label{eq:NAD_dimd}
f^{n+1}_i(\theta_k)=f^n_i(\theta_k) + \Delta t \dfrac{\mathcal F_{i+1/2}^n(\theta_k)-\mathcal F_{i-1/2}^n(\theta_k)}{\Delta w},
\ee  
for all $k=0,\dots,M$.
\begin{proposition}\label{prop:collocation_parabolic_CFL}
Under the time step restriction 
\be
\Delta t\le \dfrac{\Delta w^2}{2(U\Delta w+D)},\quad
U = \max_{i,k} |\tilde{\mathcal C}_{i+1/2}^n(\theta_k)|, 
\label{eq:nu}
\ee
the explicit scheme \eqref{eq:NAD_dimd} with flux defined by \eqref{eq:delta}-\eqref{eq:lambda_high} preserves nonnegativity for all $k=0,\dots,M$, i.e 
$ f^{n+1}_i(\theta_k)\ge 0$ if $f^n_i(\theta_k)\ge 0$, $i=0,\dots,N$, $k=0,\dots,M$. 
\end{proposition}
Higher order strong stability preserving (SSP) methods {\rm \cite{GST}} are obtained by considering a convex combination of forward Euler methods. Therefore, the non negativity result can be extended to general SSP methods. 

In practical applications, it is desirable to avoid the parabolic restriction $\Delta t = O(\Delta w^2)$ of explicit schemes. Unfortunately, fully implicit methods originate a nonlinear system of equations due to the nonlinearity of $\mathcal B[f]$ and the dependence of the weights $\delta_{i\pm 1/2}^k$ from the solution. However, we have the following nonnegativity result for the semi-implicit case
\be
f^{n+1}_i(\theta_k)=f^n_i(\theta_k) + \Delta t \dfrac{\hat{\mathcal F}_{i+1/2}^{n+1}(\theta_k)-\hat{\mathcal F}_{i-1/2}^{n+1}(\theta_k)}{\Delta w},
\label{eq:semi}
\ee
where
\be\begin{split}
\hat{\mathcal F}_{i+1/2}^{n+1}(\theta_k)=& \tilde{\mathcal C}_{i+1/2}^{k,n} \left[ (1-\delta_{i+1/2}^{k,n})f_{i+1}^{n+1}(\theta_k)+\delta^{k,n}_{i+1/2}f_i^{n+1}(\theta_k) \right]\\
&+D_{i+1/2}\dfrac{f_{i+1}^{n+1}(\theta_k)-f_i^{n+1}(\theta_k)}{\Delta w}.
\end{split}\ee
We have 
\begin{proposition}\label{prop:collocation_CFL2}
Under the time step restriction 
\be\label{eq:time_step_implicit}
\Delta t< \dfrac{\Delta w}{2U},\qquad U = \max_{i,k}|\tilde{\mathcal C}^{k,n}_{i+1/2}|
\ee
the semi-implicit scheme (\ref{eq:semi}) preserves nonnegativity, i.e 
$ f^{n+1}_i(\theta_k)\ge 0$ if $f^n_i(\theta_k)\ge 0$, $i=0,\dots,N$ for all $k=0,\dots,M$. 
\end{proposition}
We refer to \cite{PZ1} for a detailed proof. Higher order semi-implicit approximations can be constructed following \cite{BFR}. \\

\paragraph{\em Entropy property}
In order to discuss the entropy property we consider the prototype equation for all $k=0,\dots,M$
\be\label{eq:wu}
\partial_t f(\theta_k,w,t) = \partial_w \left[ P(\theta_k)(w-u)f(\theta_k,w,t) + \partial_w (D(w)f(\theta_k,w,t)) \right], 
\ee
with $w\in I = [-1,1]$ equipped with deterministic initial distribution $f(w,0)=f_0(w)$, $u=\int_{I}wf_0(w)dw\in (-1,1)$ and boundary conditions
\be\label{eq:wu_boundary}
\partial_w (D(w)f(\theta_k,w,t))+P(\theta_k)(w-u)f(\theta_k,w,t) = 0, \qquad w=\pm1.
\ee
It can be shown that the introduced structure preserving scheme dissipates the numerical entropy \cite{PZ1}

\begin{theorem}\label{th:1}
Let us consider $\mathcal B[f](\theta_k,w,t)={P(\theta_k)}(w-u)$ as in equation \eqref{eq:wu}. The numerical flux (\ref{eq:CC_flux})-(\ref{eq:f_CC}) with $\tilde{\mathcal C}_{i+1/2}^k$ and $\delta_{i+1/2}^k$ given by (\ref{eq:B_tilde})-(\ref{eq:delta}) satisfies the discrete entropy dissipation for all $k=0,\dots,M$
\be
\dfrac{d}{dt}\mathcal H_{\Delta}(f(\theta_k,w,t),f^{\infty}(\theta_k,w))=- \mathcal I_{\Delta}(f(\theta_k,w,t),f^{\infty}(\theta_k,w)),
\ee
where
\be\label{eq:relative_entropy}
\mathcal H_{\Delta w}(f(\theta_k,w,t),f^{\infty}(\theta_k,w)) = \Delta w \sum_{i=0}^N f_i \log \left(\dfrac{f_i(\theta_k,t)}{f_i^{\infty}(w,\theta_k)} \right)
\ee
and $\mathcal I_{\Delta}$ is the positive discrete dissipation function 
\be\begin{split}
& \mathcal I_{\Delta}(f(w,\theta_k,t),f^{\infty}(w,\theta_k)) = \sum_{i=0}^N
 \left[ \log \left(\dfrac{f_{i+1}(\theta_k,t)}{f^{\infty}_{i+1}(\theta_k,t)}\right)-\log\left(\dfrac{f_i(\theta_k,t)}{f_i^{\infty}(\theta_k)}\right) \right]\\
 &\qquad\qquad\cdot \left(\dfrac{f_{i+1}(\theta_k,t)}{f_{i+1}^{\infty}(\theta_k)}-\dfrac{f_i(\theta_k,t)}{f_{i}^{\infty}(\theta_k)}\right)\bar{f}_{i+1/2}^{\infty}(\theta_k)D_{i+1/2}\ge 0.
\end{split}\ee
\end{theorem}

For more general equations the above approach does not permit to prove the entropy dissipation, see \cite{PZ1}. In the following, we introduce a different class of structure preserving schemes that, in addition to preservation of the steady state of the problem, ensure the entropy dissipation.

\subsection{Entropic average schemes}
Let us consider the general class of nonlinear Fokker-Planck equation with 
gradient flow structure \cite{BCCD,CCH,CMV}
\be\label{eq:gradient_2}
\partial_t f(\theta_k,w,t) = \nabla_w \cdot [f(\theta_k,w,t)\nabla_w\xi(\theta_k,w,t)], \qquad w\in I\subseteq\RR^{d_w},
\ee
with $(\theta_k)_{k=0}^M$ the collocation nodes of the random field, and no-flux boundary conditions, where

\be\label{eq:xi_B}
\begin{split}
\nabla_w \xi(\theta_k,w,t) &= {\mathcal B}[f](\theta_k,w,t) + D {\nabla_w \log f(\theta_k,w,t)},\\
{\mathcal B}[f] (\theta_k,w,t)&=\nabla_w (U*f)(\theta_k,w,t),
\end{split}\ee
with $U(\theta_k,\cdot)$ an uncertain interaction potential. A stochastic free energy functional is defined as follows
\[
\begin{split}
\mathcal E(\theta_k,t) = \dfrac{1}{2}\int_{\RR^d}(U*f)(\theta_k,w,t)f(\theta_k,w,t)dw+D\int_{\RR^d}\log f(\theta_k,w,t) f(\theta_k,w,t)dw.
\end{split}
\]
which is dissipated along solutions as
\be\label{eq:dissipation}
\dfrac{d}{dt}\mathcal E(\theta_k,t) = -\int_{\RR^d}|\nabla_w\xi|^2 f(\theta_k,w,t)dw =- \mathcal I(\theta_k,t),
\ee
where $\mathcal I(\theta_k,\cdot)$ is the entropy dissipation function. 
The corresponding discrete free energy is given by
\be\begin{split}\label{eq:discrete_entropy}
&\mathcal E_{\Delta}(\theta_k,t) = \Delta w \sum_{j=0}^N \Big[ \dfrac{1}{2}\Delta w\sum_{i=0}^N U_{j-i}(\theta_k)f_i(\theta_k,t)f_j(\theta_k,t)\\
&\qquad\qquad +D f_j(\theta_k,t)\log f_j(\theta_k,t) \Big]
\end{split}\ee

In this case it is not possible to show that the discrete entropy functional \eqref{eq:discrete_entropy} is dissipated by the SP--CC type schemes developed in the previous sections, see \cite{PZ1}. For this reason we introduce the new entropic family of flux function
\be\label{eq:entropic_average}
 \tilde f^{E}_{i+1/2}(\theta_k,t) =
 \begin{cases}
  \dfrac{f_{i+1}(\theta_k,t)-f_i(\theta_k,t)}{\log f_{i+1}(\theta_k,t)-\log f_i(\theta_k,t)} & f_{i+1}(\theta_k,t)\ne f_i(\theta_k,t),\\
 f_{i+1}(\theta_k,t) & f_{i+1}(\theta_k,t)=f_i(\theta_k,t),
 \end{cases}\ee
for all $k=0,\dots,M$.  We will refer to the above approximation of the solution at the grid point $i+1/2$ as \emph{entropic average} of the grid points $i$ and $i+1$. In the general case of the flux function \eqref{eq:flux} with non constant diffusion the resulting numerical flux reads
\be\begin{split}\label{eq:F_E}
&\mathcal F^E_{i+1/2}(\theta_k,t) = D_{i+1/2}\Bigg( \frac{\tilde {\mathcal C}_{i+1/2}(\theta_k,t)}{D_{i+1/2}}\\
&\qquad\qquad\qquad\qquad+\dfrac{\log f_{i+1}(\theta_k,t)-\log f_i(\theta_k,t)}{\Delta w}  \Bigg) \tilde f^E_{i+1/2}(\theta_k,t).
 \end{split}\ee
Concerning the stationary state, we obtain immediately by imposing the numerical flux equal to zero
\[
\frac{\tilde {\mathcal C}_{i+1/2}(\theta_k,t)}{D_{i+1/2}}+\dfrac{\log f_{i+1}(\theta_k,t)-\log f_i(\theta_k,t)}{\Delta w}=0,
\]  
and therefore we get
\be
\frac{f_{i+1}(\theta_k,t)}{f_i(\theta_k,t)} = \exp\left(-\frac{\Delta w\,\tilde {\mathcal C}_{i+1/2}(\theta_k,t)}{D_{i+1/2}}\right).
\ee 
By equating the above ratio with the quasi-stationary approximation \eqref{eq:quasi_SS} we get the same expression for $\tilde{\mathcal{C}}_{i+1/2}(\theta_k,t)$ for all $k=0,\dots,M$ as in \eqref{eq:B_tilde}
\be\label{eq:B_tilde2}
\tilde{\mathcal{C}}_{i+1/2}(\theta_k,t)=\dfrac{D_{i+1/2}}{\Delta w}\int_{w_i}^{w_{i+1}}\dfrac{\mathcal B[f](w,\theta_k,t)+D'(w)}{D(w)}dw.
\ee 
A fundamental result concerning the entropic average \eqref{eq:entropic_average} is the following 
 \begin{lemma}
The entropy average defined in \eqref{eq:entropic_average} may be written as a convex combination with nonlinear weights
\be
\tilde{f}^{E}_{i+1/2}(\theta_k,t) = \delta_{i+1/2}^{k,E} f_i(\theta_k,t) +(1-\delta_{i+1/2}^{k,E})f_{i+1}(\theta_k,t),
\ee
where 
\be\label{eq:deltaE}
\delta_{i+1/2}^{k,E} = \dfrac{f_{i+1}(\theta_k,t)}{f_{i+1}(\theta_k,t)-f_i(\theta_k,t)}+\dfrac{1}{\log f_i(\theta_k,t)-\log f_{i+1}(\theta_k,t)}\in(0,1).
\ee
 \end{lemma}
 \begin{remark} As a consequence the Chang-Cooper type average \eqref{eq:f_CC} and the entropic average \eqref{eq:entropic_average} define the same quantity at the steady state when $f_i(\theta_k,t)=f_i^{\infty}(\theta_k)$. 
 \end{remark}
On the contrary to the Chang-Cooper average the restrictions for the non negativity property of the solution are stronger.
Therefore, similar to central differences, we have a restriction on the mesh size which becomes prohibitive for small values of the diffusion function $D(w)$. It is possible to show that the same condition is necessary also for the non negativity of semi-implicit approximations.

Concerning the entropy dissipation we can summarize the main results in the following \cite{PZ1}
\begin{theorem}\label{pr:2}
The numerical flux \eqref{eq:F_E}-\eqref{eq:entropic_average} for a constant diffusion $D$ satisfies the discrete entropy dissipation
\be
\dfrac{d}{dt}\mathcal{E}_{\Delta}(\theta_k,t)=- \mathcal I_{\Delta}(\theta_k,t),
\ee
where $\mathcal{E}_{\Delta}(\theta_k,t)$ is given by \eqref{eq:discrete_entropy} 
and $I_{\Delta}(\theta_k,t)$ is the discrete entropy dissipation function 
\be
\mathcal I_{\Delta}(\theta_k,t) = \Delta w \sum_{j=0}^N  (\xi_{j+1}(\theta_k,t)-\xi_j(\theta_k,t))^2  \tilde{f}^{E}_{i+1/2}(\theta_k,t) \ge 0,
\ee
with $\xi_{j+1}(\theta_k,t)-\xi_j(\theta_k,t)$ the discrete version of \eqref{eq:xi_B}.
\end{theorem}
Further, we can state the following entropy dissipation results for problem \eqref{eq:wu} in the nonlogarithmic Landau form \eqref{eq:Landau_nonlog}.
\begin{theorem}\label{th:1b}
Let us consider $\mathcal B[f](w,\theta_k,t)=P(\theta_k)(w-u)$ as in equation \eqref{eq:wu}. The numerical flux \eqref{eq:F_E}-\eqref{eq:entropic_average} with $\tilde{\mathcal C}_{i+1/2}^k$ given by (\ref{eq:B_tilde}) satisfies the discrete entropy dissipation
\be
\dfrac{d}{dt} \mathcal H_{\Delta}(f(\theta_k,t),f^{\infty}(\theta_k,t))  = -\mathcal I^E_{\Delta}(f(\theta_k,t),f^{\infty}(\theta_k,t)),
\ee
where $\mathcal H_{\Delta w}(f(\theta_k,t),f^{\infty}(\theta_k,t))$ is given by \eqref{eq:relative_entropy} and $\mathcal I^E_{\Delta}(\theta_k,t)$ is the positive discrete dissipation function 
\be\begin{split}
\mathcal I^E_{\Delta}(f(\theta_k,t),f^{\infty}(\theta_k,t))=& \sum_{i=0}^N \left[ \log \left(\dfrac{f_{i+1}(\theta_k,t)}{f^{\infty}_{i+1}(\theta_k,t)}\right) -\log \left(\dfrac{f_i(\theta_k,t)}{f^{\infty}_i(\theta_k,t)}\right) \right]^2\\
&\cdot D_{i+1/2} \tilde f^E_{i+1/2}(\theta_k,t)\ge 0.
\end{split}\ee
\end{theorem}

\subsection{Numerical results}\label{sec:FP_collocation}
We report a numerical example obtained with the collocation approach in combination with the structure--preserving numerical methods. We consider a stochastic Fokker--Planck equation with uncertainty in the initial distribution, i.e. 
\be\begin{cases}\label{eq:FP_num_SP}
\partial_t f(\theta,w,t) = \partial_w \Big[ wf(\theta,w,t)+{T(\theta)}\partial_w^2 f(\theta,w,t) \Big], \\
f(\theta,w,0) = f_0(\theta,w),
\end{cases}\ee
for all $w\in\RR$ with 
\be\label{eq:FP_initial}
f_0(\theta,w) = \dfrac{1}{2} \Bigg\{ \dfrac{1}{\sqrt{2\pi\sigma^2(\theta)}}e^{-\frac{(w-c)^2}{2\sigma^2(\theta)}}+\dfrac{1}{\sqrt{2\pi\sigma^2(\theta)}}e^{-\frac{(w+c)^2}{2\sigma^2(\theta)}} \Bigg\}, \qquad c = 1/10
\ee
and $\sigma^2(\theta)=1/10+\epsilon\theta$, $\theta\sim U([-1,1])$, $\epsilon =5\times 10^{-3}$. In \eqref{eq:FP_num_SP} the diffusion coefficient is the temperature
\[
T(\theta) = \int_{\RR} w^2f_0(\theta,w)dw.
\]
It is well--known that the steady--state solution of this problem is the Maxwellian distribution \eqref{eq:maxw}.   

\begin{figure}[t]
\centering
\includegraphics[scale=0.5]{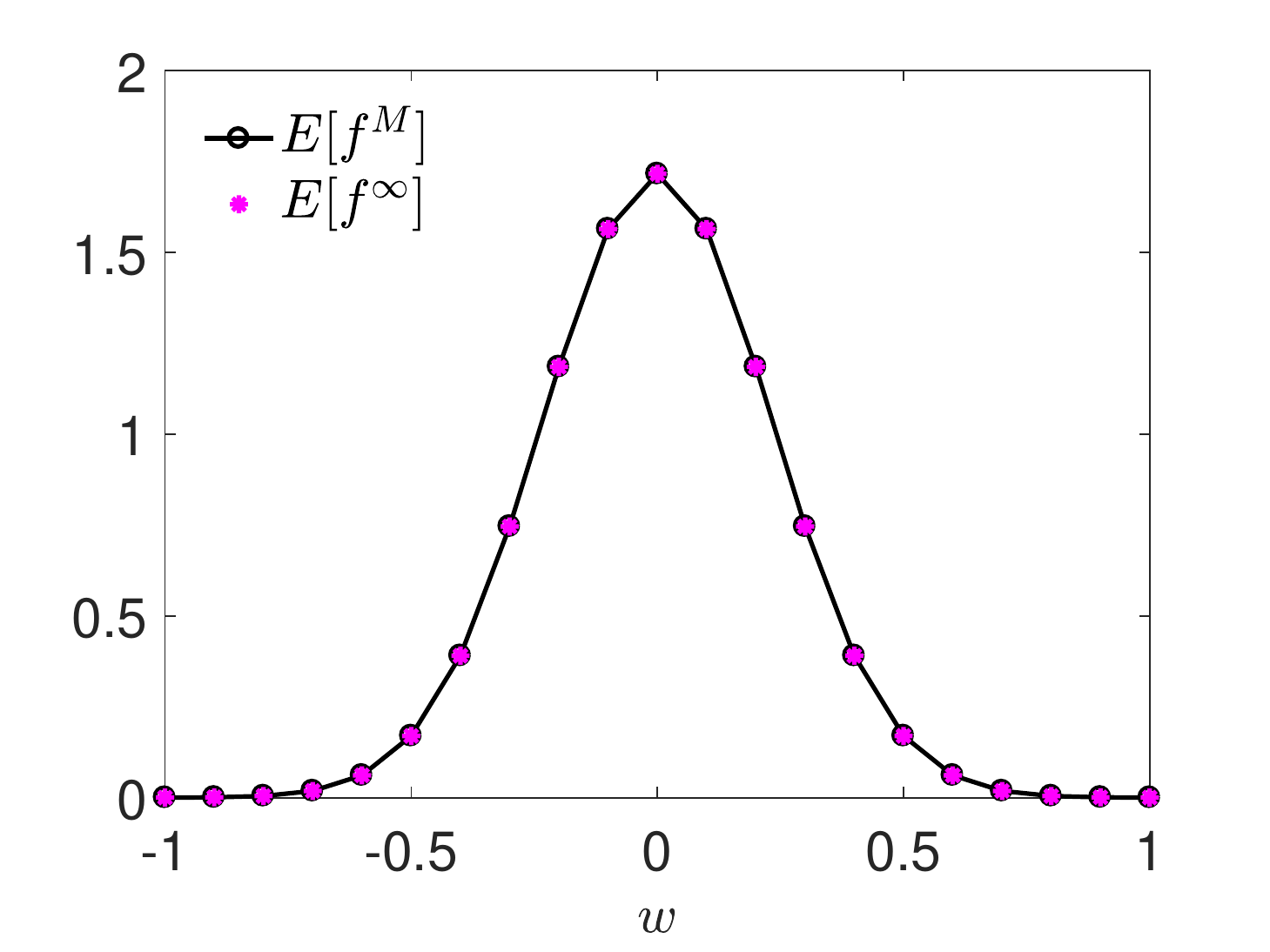}
\includegraphics[scale=0.5]{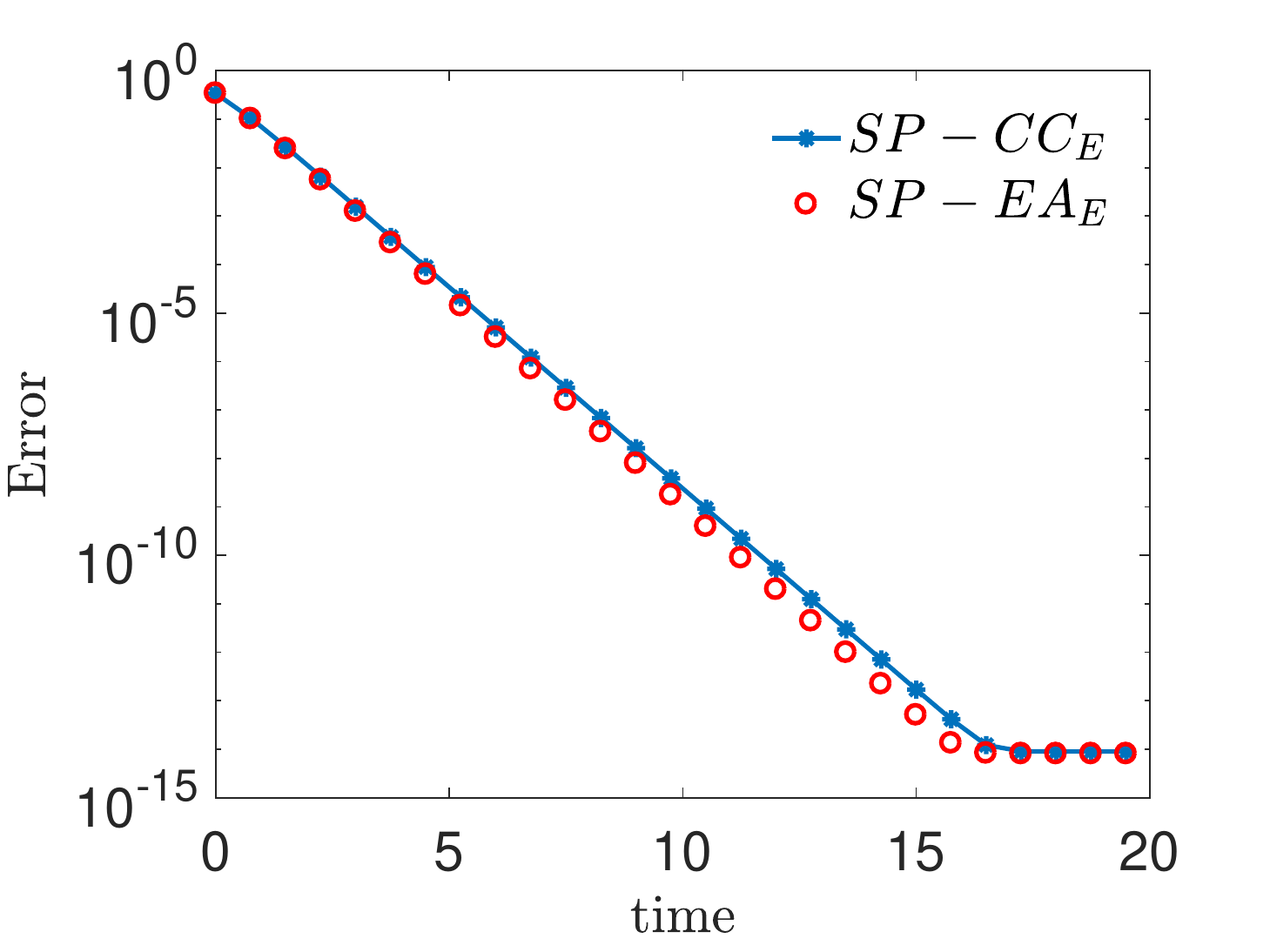}
\caption{Left: exact and numerical approximation of the expected steady state distribution. Right: evolution of the $L^1$ relative error for the expected solution calculated for both $SP-CC_E$ and $SP-EA_E$ methods. In both figures we considered a grid on $[-1,1]$ with $N=21$ points and $M=10$ nodes in the random field, the final time $T=20$ and $\Delta t=\Delta w^2/2$. The nodes of the random field have been chosen with Gauss--Legendre polynomials.  }
\label{fg:fig1}
\end{figure}

\begin{figure}[h]
\centering
\includegraphics[scale=0.5]{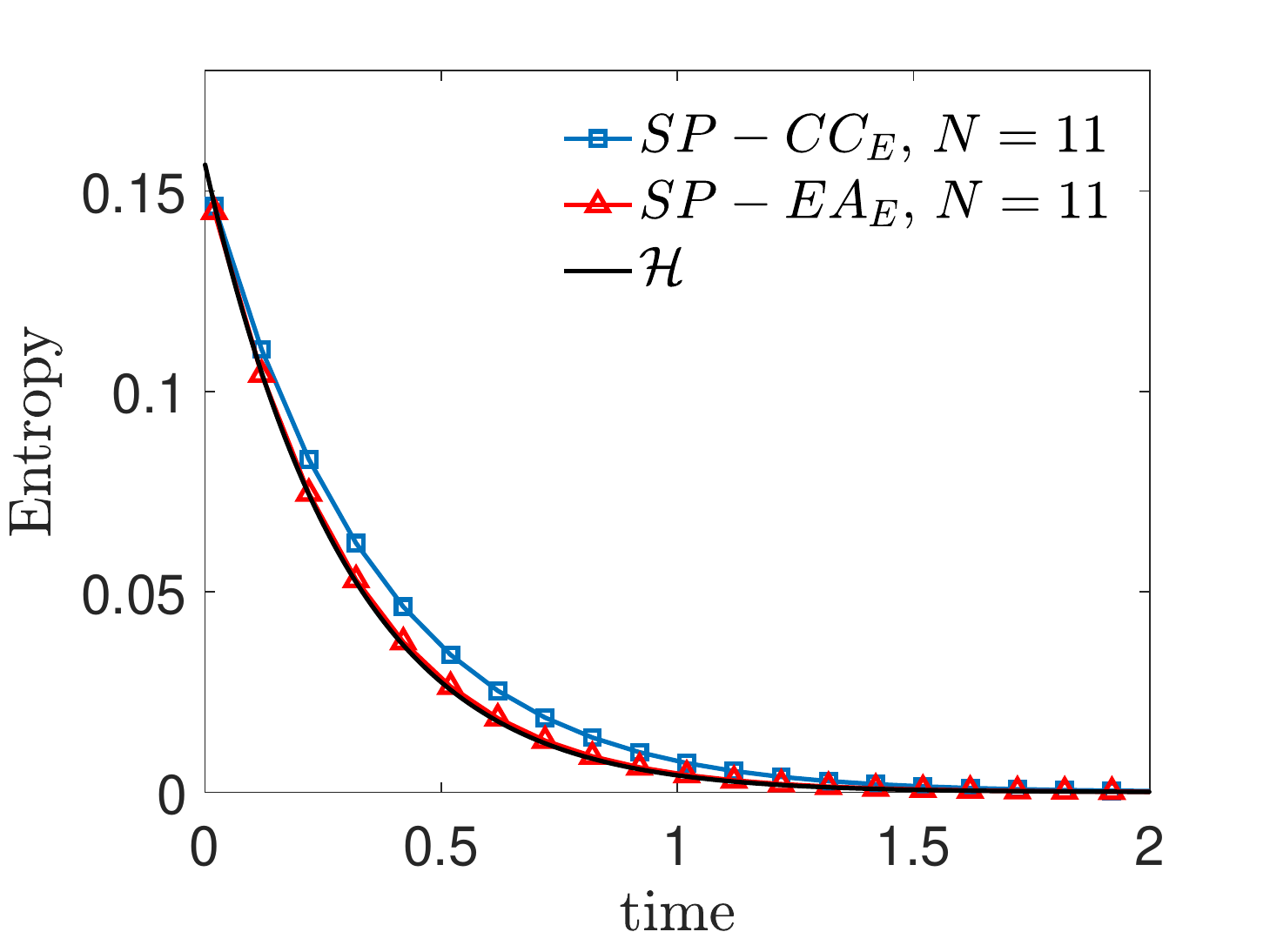}
\caption{Dissipation of the numerical expected entropy for $SP-CC_E$ and $SP-EA_E$ schemes on a coarse grid with $N=11$ points.}
\label{fg:fig2}
\end{figure}

In the previous paragraphs we showed how an essential aspect for the accurate description of the stochastic steady state relies in the approximation of the family of integrals $\delta_{i+1/2}^k$, $\lambda_{i+1/2}^k$, see \eqref{eq:delta}--\eqref{eq:lambda_high}. In this case, however, since the steady state is known we can evaluate exactly these weight functions as in (\ref{eq:lambda_inf})-(\ref{eq:delta_inf}). We postpone to the last section of the present contribution the discussion on numerical results obtained with more general weight functions for which no exact formulation are given. In Figure \ref{fg:fig1} (right) we report the relative $L^1$ error for the expectation of the solution in time. As expected the schemes are capable to capture the stochastic steady state exactly. Next in Figure \ref{fg:fig2} the evolution of the expectation for the numerical entropy is given.

\section{Variance reduction Monte Carlo methods}\label{sec:4}

Among the different type of techniques used in UQ, certainly the Monte Carlo methods
represent one of the most popular and important class \cite{Caflisch, DP, Giles, PT2}. They show all their
potential when the dimension of the uncertainty space becomes very large. In addition, Monte Carlo methods are effective when the probability distribution of the random inputs is not known analytically or lacks of regularity since other approaches based on orthogonal stochastic polynomials may be impossible to use of may produce poor results.     

In this section, we first describe a standard Monte Carlo approach which deals with random initial data and then we describe a modification of this algorithm which permits to strongly decrease the computational costs and increase the accuracy close to the equilibrium steady state. 

\subsection{The standard Monte Carlo method}
 We describe the method when applied to the solution of a Vlasov-Fokker-Planck type equation of the type (\ref{eq:MF_general}) with deterministic parameters $P=P(x,x_*,w,w_*)$ and $D=D(w,t)$ and random initial data $f(\theta,x,w,0)=f_0(\theta,v,w)$.
First we assume that the kinetic equation has been discretized by a deterministic solver in the variables $w$, $x$ and $t$. In this setting, the simplest Monte Carlo (MC) method for UQ in kinetic equation is based on the following steps. 

\begin{algorithm}[Standard Monte Carlo (MC)  method]~
\begin{enumerate}
\item {\bf Sampling}: Sample $M$ independent identically distributed (i.i.d.) initial data $f_0^k$, $k=1,\ldots,M$ from
the random field $f_0$ and approximate these over the grid (for example by piece-wise constant cell averages).
\item {\bf Solving}: For each realization $f_0^k$ the underlying kinetic equation (\ref{eq:MF_general}) is solved numerically by the deterministic solver. We denote the solutions at time $t^n$ by $f^{k,n}_{\Delta w,\Delta x}$, $k=1,\ldots,M$, where $\Delta w$ and $\Delta x$ characterizes the discretizations  in $w$ and $x$. 
\item {\bf Estimating}: Estimate the expected value of the random solution field with the sample mean of the approximate solution
\be
E_M[f^{n}_{\Delta w}]=\frac1{M} \sum_{k=1}^M f^{k,n}_{\Delta w, \Delta x}.
\label{mcest}
\ee
\end{enumerate}
\end{algorithm}
The above algorithm is straightforward to implement in any existing code for the Vlasov-Fokker-Planck equations. Furthermore, the only (data) interaction between different samples { is in step $3$}, when ensemble averages are computed. Thus, the MC algorithms for UQ are non-intrusive {and easily parallelizable as well}.

The typical error estimate that one obtains using such an approach is of the type
\be
\|E[f(\cdot,t^n)]-E_M[f^{n}_{\Delta w}]\| \leq C_1 M^{-1/2} + C_2 (\Delta w)^q + C_3 (\Delta x)^p + C_4 (\Delta t)^r
\ee
where $\|\cdot\|$ is a suitable norm, $C_1$, $C_2$, $C_3$ and $C_4$ are positive constants depending only on the
second moments of the initial data and the interaction term, and $q$, $p$ and $r$ characterize the accuracy of the discretizations in the phase-space. Clearly, it is possible to equilibrate the discretization and the sampling errors in the a-priori estimate taking $M=O(\Delta w^{-2q})$,
$\Delta x=O(\Delta w^{q/p})$ and $\Delta t=O(\Delta w^{q/r})$. This means that in order to have comparable errors the number of samples should be extremely large, especially when dealing with high order deterministic discretizations. This may make the collocation Monte Carlo approach very expensive in practical applications.  
In order to address the slow convergence of Monte Carlo methods, we discuss in the next paragrph the development of variance reduction Monte Carlo methods.

\subsection{The Micro-Macro Monte Carlo method}
In order to improve the performances of standard MC methods, we introduce a novel class of variance reduction Monte Carlo methods \cite{DP}. The key idea on which they rely is to take advantage of the knowledge of the steady state solution of the kinetic equation in order to reduce both the variance and the computational cost of the Monte Carlo estimate. The method is essentially a control variate strategy based on a suitable microscopic-macroscopic decomposition of the distribution function. 

We describe the method in the space homogeneous case, an example of such technique in the non homogeneous case is reported in Section \ref{sec:6}, while we refer to \cite{DP} for a detailed discussion and extensions of such method to more general kinetic equations. Following Section \ref{sec:micro-macro} we introduce the  Micro--Macro decomposition $$f(\theta,w,t)=f^{\infty}(\theta,w)+g(\theta,w,t),$$ where $f^{\infty}(\theta,w)$ is the steady state solution of the problem considered. Then, the idea consists in using the Monte Carlo estimation procedure only to the non equilibrium part $g(\theta,w,t)$ solution of (\ref{eq:micro-macro}).

The crucial aspect is that the equilibrium state $g^\infty(\theta,w)$ is zero and therefore, independent from $\theta$. More precisely, we can decompose the expected value of the distribution function in
an equilibrium and non equilibrium part 
\be
\begin{split}
{\mathbb E}[f](w,t)]&=\int_{\Omega} f(\theta,w,t)p(\theta)d\theta\\
&=\int_{\Omega} f^\infty(\theta,w)p(\theta)d\theta+
\int_{\Omega} g(\theta,w,t)p(\theta)d\theta,
\end{split}\ee
and then exploit the fact that $f^\infty(\theta,w)$ is known to have an estimate of the error committed by the Monte Carlo integration of 
type
\be 
e_M[f]\simeq \sigma_{g}M^{-1/2}
\ee
instead of 
\be 
e_M[f]\simeq \sigma_{f}M^{-1/2},
\ee
where $\sigma_g$ and $\sigma_f$ are the variances of respectively the perturbation and the distribution function and where we have supposed for simplicity that expected value of the equilibrium part is computed with a negligible error. Now, since it is known that the perturbation $g$ goes to zero in time exponentially fast, then also its variance goes to zero, which means that in the steady state solution limit the Monte Carlo integration becomes only dependent on the way in which the expected value of the equilibrium part is computed.

The simplest version of the algorithm consists of the following steps:
\begin{algorithm}[Micro-Macro Monte Carlo (M$^3$C) method]~
\begin{enumerate}
\item {\bf Small scale sampling}: Sample $M_E$ independent identically distributed (i.i.d.) initial data $f_0^k$, $k=1,\ldots,M_E$ from
the random field $f_0$. For each sample compute the corresponding equilibrium state $f_{\Delta w}^{k,\infty}$ from its moments evaluated through suitable quadrature rules in $w$ based on the discretization parameter $\Delta w$.
\item {\bf Large scale sampling}: Select $M \ll M_E$ samples $f_0^k$, $k=1,\ldots,M$ and compute $g_0^k=f_0^k-f^{k,\infty}$ and approximate
these over the grid (for example by piece-wise constant cell averages). 
\item {\bf Solving}: For each realization $g_0^k$ the underlying kinetic equation (\ref{eq:micro-macro}) is solved numerically by the deterministic solver. We denote the solutions at time $t^n$ by $g^{k,n}_{\Delta w}$, $k=1,\ldots,M$.

\item {\bf Estimating}: We estimate the expected value of the random solution field 
\[
f^{k,n}_{\Delta w}=f_{\Delta w}^{k,\infty}+g^{k,n}_{\Delta w},
\]
with the sample mean of the approximate solution
\be
E_{M,M_E}[f^{n}_{\Delta w}]=\frac1{M_E} \sum_{k=1}^{M_E} f_{\Delta w}^{k,\infty}+\frac1{M} \sum_{k=1}^M g^{k,n}_{\Delta w}.
\label{mcest2}
\ee
\end{enumerate}
\end{algorithm}
Using such an approach one obtains an error estimate of the type
\be
\|E[f(\cdot,t^n)]-E_{M,M_E}[f^{n}_{\Delta w}]\| \leq C_E M_E^{-1/2}+C^n_1 M^{-1/2} + C_2 (\Delta w)^q +C_3(\Delta t)^r
\ee
where now the constant $C^n_1$ depends on time and on the second moment of the solution $g(\theta,w,t^n)$ which will  vanish for large times. In fact, independently of $\theta$ we have that $g(\theta,w,t^n)\to 0$ as $n\to \infty$. 
Therefore, the method reduces the variance of the estimator in time and asymptotically, since $C^n_1\to 0$ as $n\to \infty$, depends only on the fine scale sampling which does not affect the overall computational cost.

The efficiency of the M$^3$C can be further improved in the case of monotonic convergence to equilibrium of the distribution function $f$ by introducing a strategy of sampling reduction at each time step. The resulting algorithm is the following  

\begin{algorithm}[Fast Micro-Macro Monte Carlo (FM$^3$C) method]~
\begin{enumerate}
\item {\bf Small scale sampling}: Sample $M_E$ independent identically distributed (i.i.d.) initial data $f_0^k$, $k=1,\ldots,M_E$ from
the random field $f_0$. For each sample compute the corresponding equilibrium state $f_{\Delta w}^{k,\infty}$ from its moments evaluated through suitable quadrature rules in $w$ based on the discretization parameter $\Delta w$.
\item {\bf Large scale sampling}: Select $M_0 \ll M_E$ samples $f_0^k$, $k=1,\ldots,M_0$ and compute $g_0^k=f_0^k-f^{k,\infty}$ and approximate
these over the grid (for example by piece-wise constant cell averages). 
\item {\bf Solving}: For each realization $g_0^k$ the underlying kinetic equation (\ref{eq:micro-macro}) is solved numerically by the deterministic solver. This is realized at each time step $n=0,1,2,\ldots$ as follows.
\begin{enumerate}
\item {\bf Advance in time}: Starting from $g^{k,{n}}_{\Delta w}$, $k=1,\ldots,M_n$ compute the solution $g^{k,{n+1}}_{\Delta w}$ with one time step of the deterministic solver.
\item {\bf Discard samples}: At each time step we compute the variance of $g^{k,{n+1}}_{\Delta w}$ as
\[
{\it Var}_{M_{n}}[g^{n+1}_{\Delta w}]=\frac1{M_{n}} \sum_{k=1}^{M_n} (g^{k,{n+1}}_{\Delta w}-E_{M_{n}}[g^{k,n+1}_{\Delta w}])^2 \leq {\it Var}_{M_n}[g^{n}_{\Delta w}].
\]
Set $M_{n+1}=[\![M_n \left({\it Var}_{M_{n}}[g^{n+1}_{\Delta w}]/{\it Var}_{M_{n}}[g^{n}_{\Delta w}]\right)]\!]$ where $[\![\cdot]\!]$ denotes the integer part and discard uniformly $M_{n}-M_{n+1}$ samples.
\end{enumerate} 
\item {\bf Estimating}: We estimate the expected value of the random solution field 
\[
f^{k,n}_{\Delta w}=f_{\Delta w}^{k,\infty}+g^{k,n}_{\Delta w},
\]
with the sample mean of the approximate solution
\be
E_{M_n,M_E}[f^{n}_{\Delta w}]=\frac1{M_E} \sum_{k=1}^{M_E} f_{\Delta w}^{k,\infty}+\frac1{M_n} \sum_{k=1}^{M_n} g^{k,n}_{\Delta w}.
\label{mcest2fast}
\ee
\end{enumerate}
\end{algorithm}
The algorithm preserves the advantages of the simple M$^3$C method but with a greater computational efficiency since the number of samples, and therefore the number of deterministic equations that we have to solve, decreases in time and asymptotically vanishes.

\begin{remark}~
\begin{itemize}
\item In the case the underlying uncertainty probability density function $p(\theta)$ is known, the M$^3$C method can be applied without any small scale sampling since the estimate of the expected value reduces to
\be
E_{M_n}[f^{n}_{\Delta w}]=\int_{\Omega} f_{\Delta w}^{\infty}(\theta,w)p(\theta)d\theta+\frac1{M_n} \sum_{k=1}^{M_n} g^{k,n}_{\Delta w}.
\label{mcest2b}
\ee  
In this case M$^3$C methods achieves arbitrary accuracy for large times.
\item In contrast with Multi Level Monte Carlo (MLMC) methods \cite{Giles}, which can produce non monotone estimators, the estimators produced by the M$^3$C method are monotonic, i.e. mean estimator of positive quantities (such as density) is also positive, the same holds true for the entropy property.
\item The extension of the M$^3$C method to the non homogeneous case is straightforward, whereas the advantages of the FM$^3$C method depend on the type of problem considered. Applications are reported in Section \ref{sec:6}. For more general cases we refer to \cite{DP} where e detailed discussion is done.
\end{itemize}
\end{remark}

\subsection{Numerical results}
In this section we show some results concerning the Micro-Macro Monte Carlo methods by comparing them to the standard MC method for UQ. In particular, we study the behaviors of our approach in solving the stochastic Fokker--Planck equation with uncertainty in the initial distribution (\ref{eq:FP_num_SP})-(\ref{eq:FP_initial}).
\begin{figure}[t]
\centering
\includegraphics[scale=0.5]{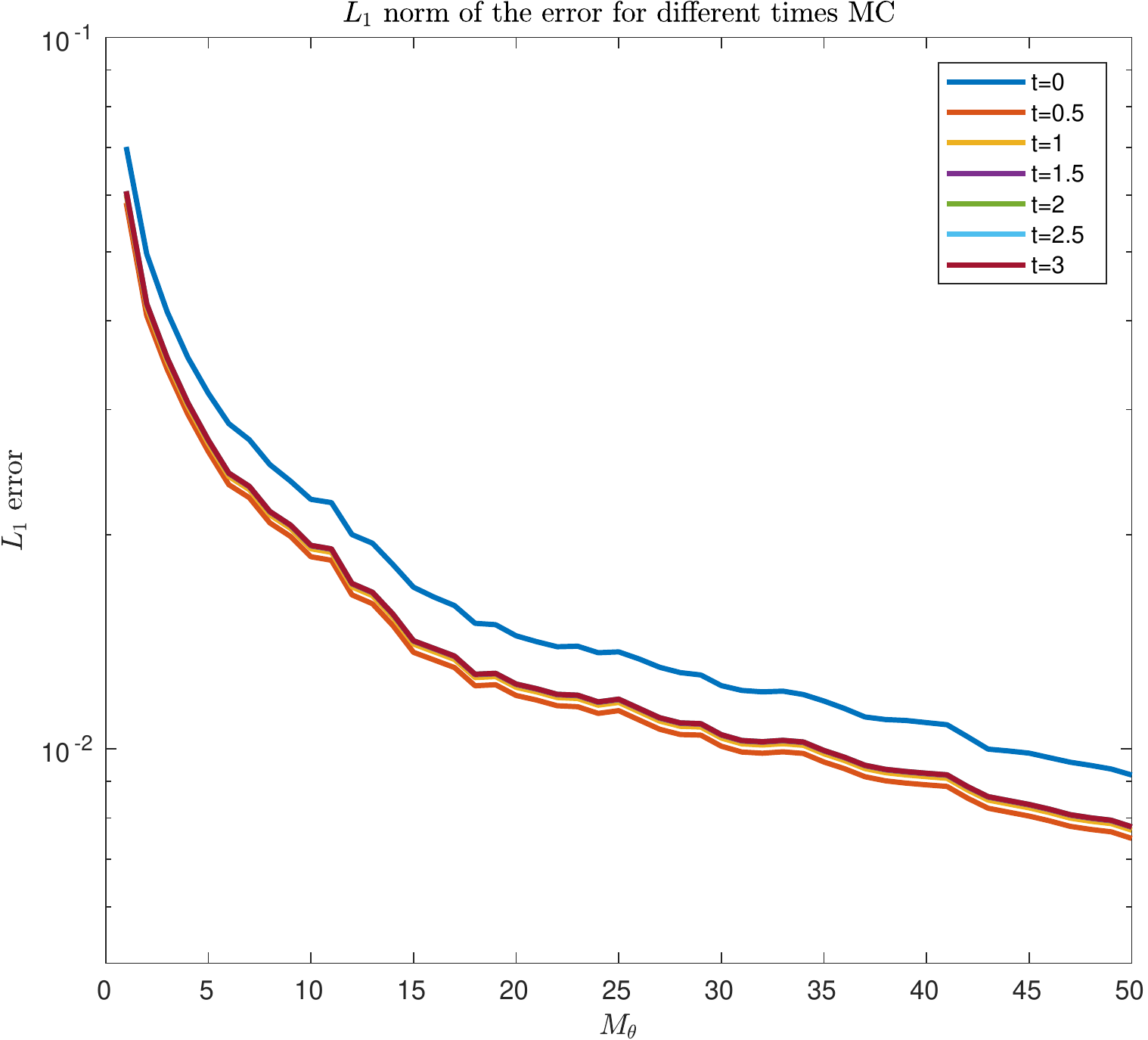}
\includegraphics[scale=0.5]{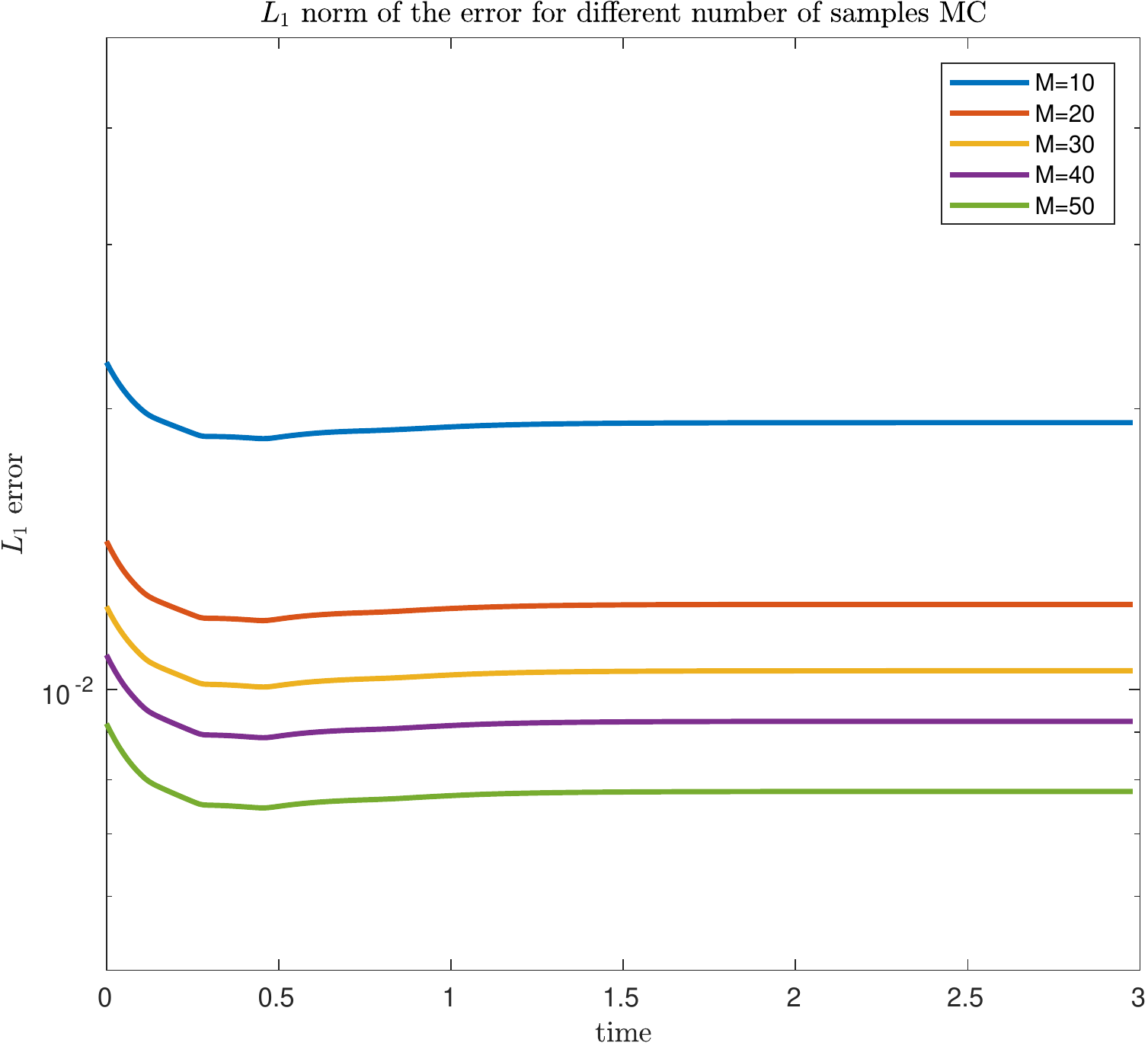}
\caption{Monte Carlo method. Left: evolution of the $L_1$ norm of the error for the expected distribution over an increasing number of stochastic inputs. Different lines represent the error at different times. Right: evolution of the $L^1$ norm of the error for the expected distribution over time. Different lines represent a different number of stochastic inputs over which the expected value is computed. Grid $[-1,1]$ with $N=100$ points, final time $T=3$ and $\Delta t=\Delta w^2/2$. The solution has been averaged over 100 different realizations.}
\label{fig1Mc}
\end{figure}

\begin{figure}[t]
\centering
\includegraphics[scale=0.5]{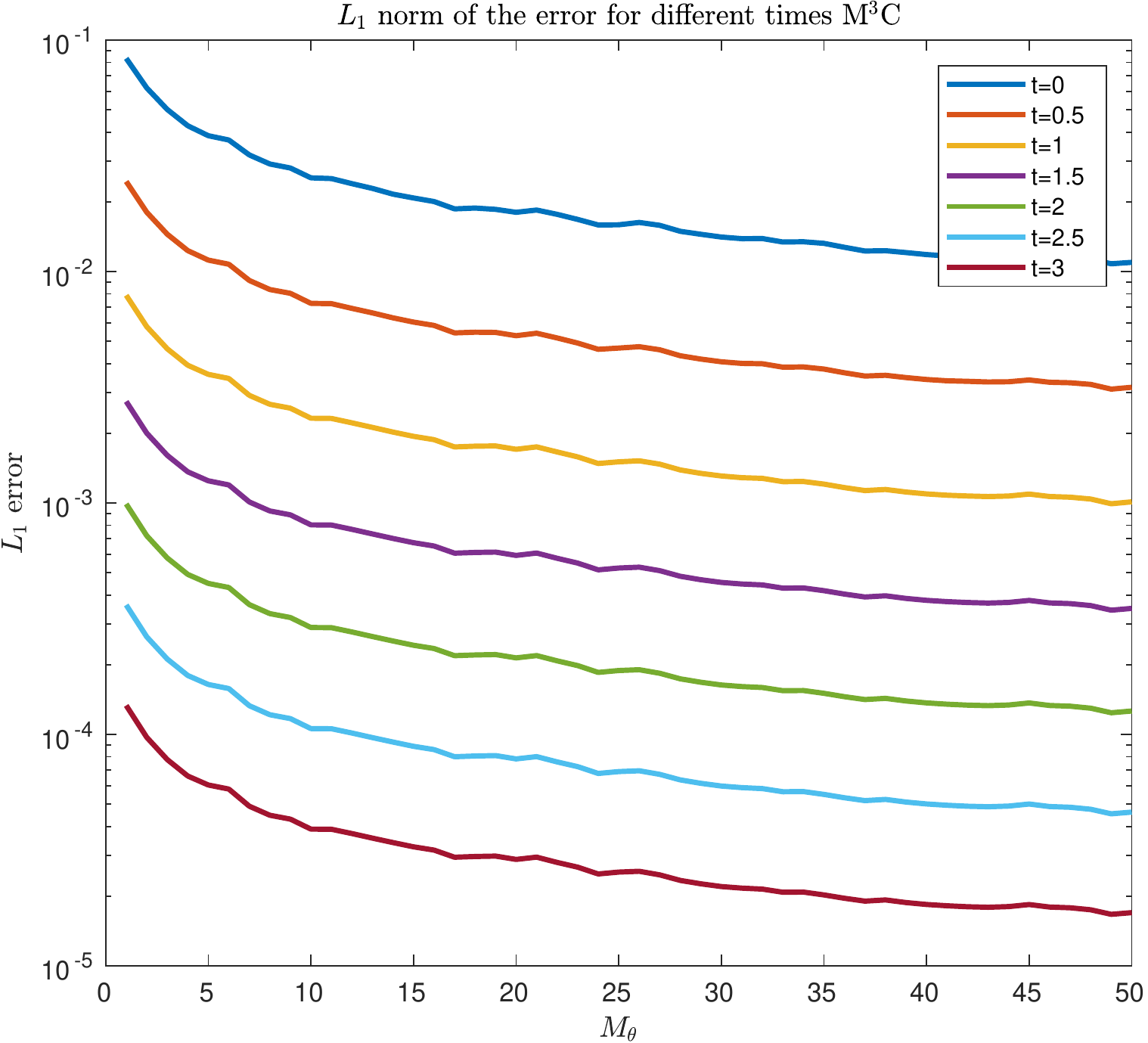}
\includegraphics[scale=0.5]{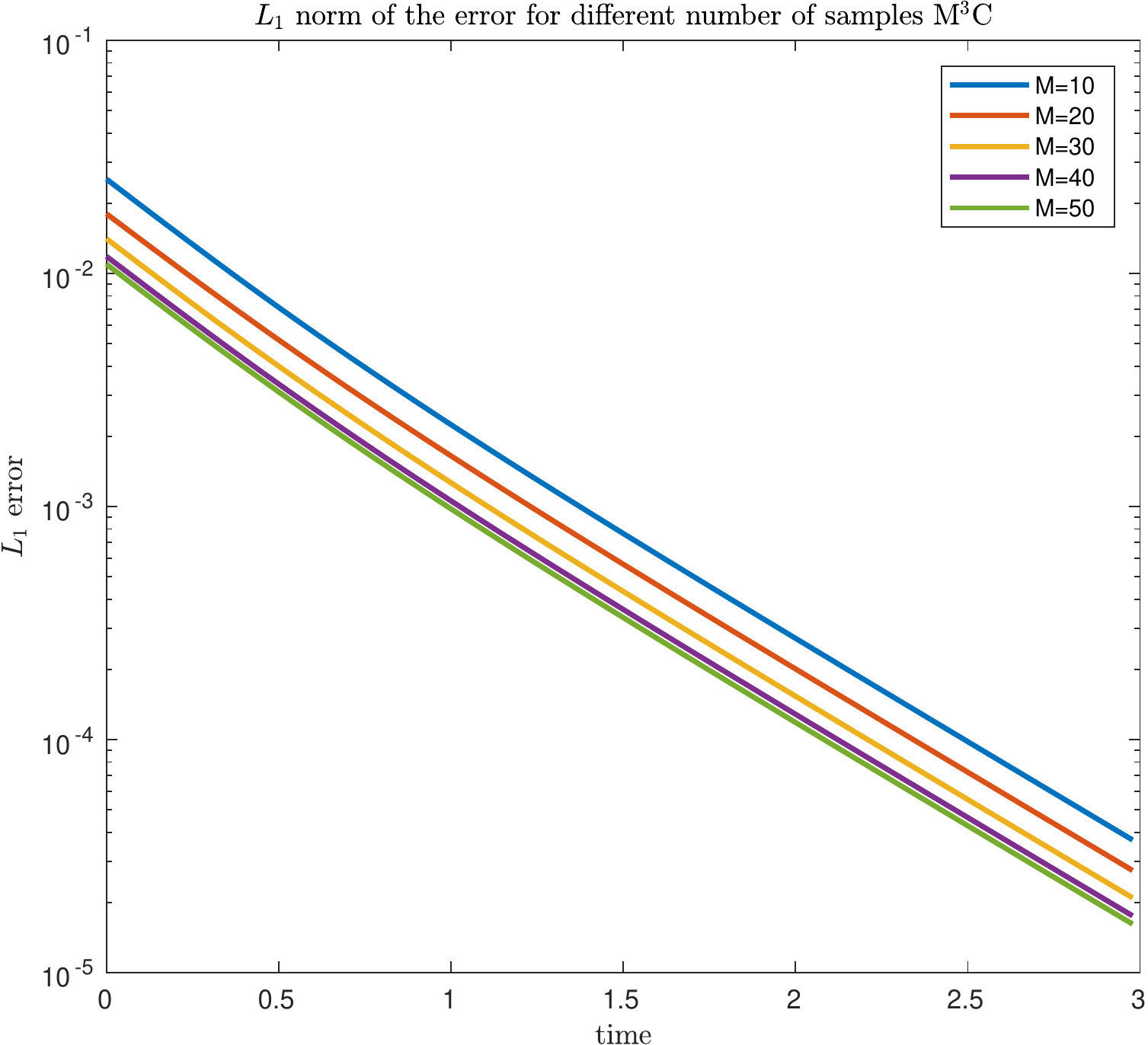}
\caption{M$^3$C method. Left: evolution of the $L_1$ norm of the error for the expected distribution over an increasing number of stochastic inputs. Different lines represent the error at different times. Right: evolution of the $L^1$ norm of the error for the expected distribution over time. Different lines represent a different number of stochastic inputs over which the expected value is computed. Grid $[-1,1]$ with $N=100$ points, final time $T=3$ and $dt=dw^2/2$. The solution has been averaged over 100 different realizations.}
\label{fig2Mc}
\end{figure}

In Figure \ref{fig1Mc} we report the $L_1$ norm of the error for the expected solution with a standard MC method. Left image shows the error for an increasing number of samples for different times, while right image shows the trend of the error over time for a different number of random inputs. The final time is set to $T_f=3$, the number of cells in velocity is $100$, while the stability condition gives $\Delta t=0.5\Delta w^2$. The maximum number of samples which furnishes the set of initial conditions is $M_\theta=50$, while the solution is averaged over 100 realizations. One can clearly see the $M_\theta^{-1/2}$ slope for the error in the left picture. 

In Figure \ref{fig2Mc}, the $L_1$ norm of the error is reported in the same setting of the MC case for the M$^3$C method both for the left and the right images. The same number of averages and stochastic initial condition is employed. We can see how the error decreases as a function of time in an exponential fashion at the contrary of the MC case for which the error is almost independent on time. 

\begin{figure}[t]
\centering
\includegraphics[scale=0.5]{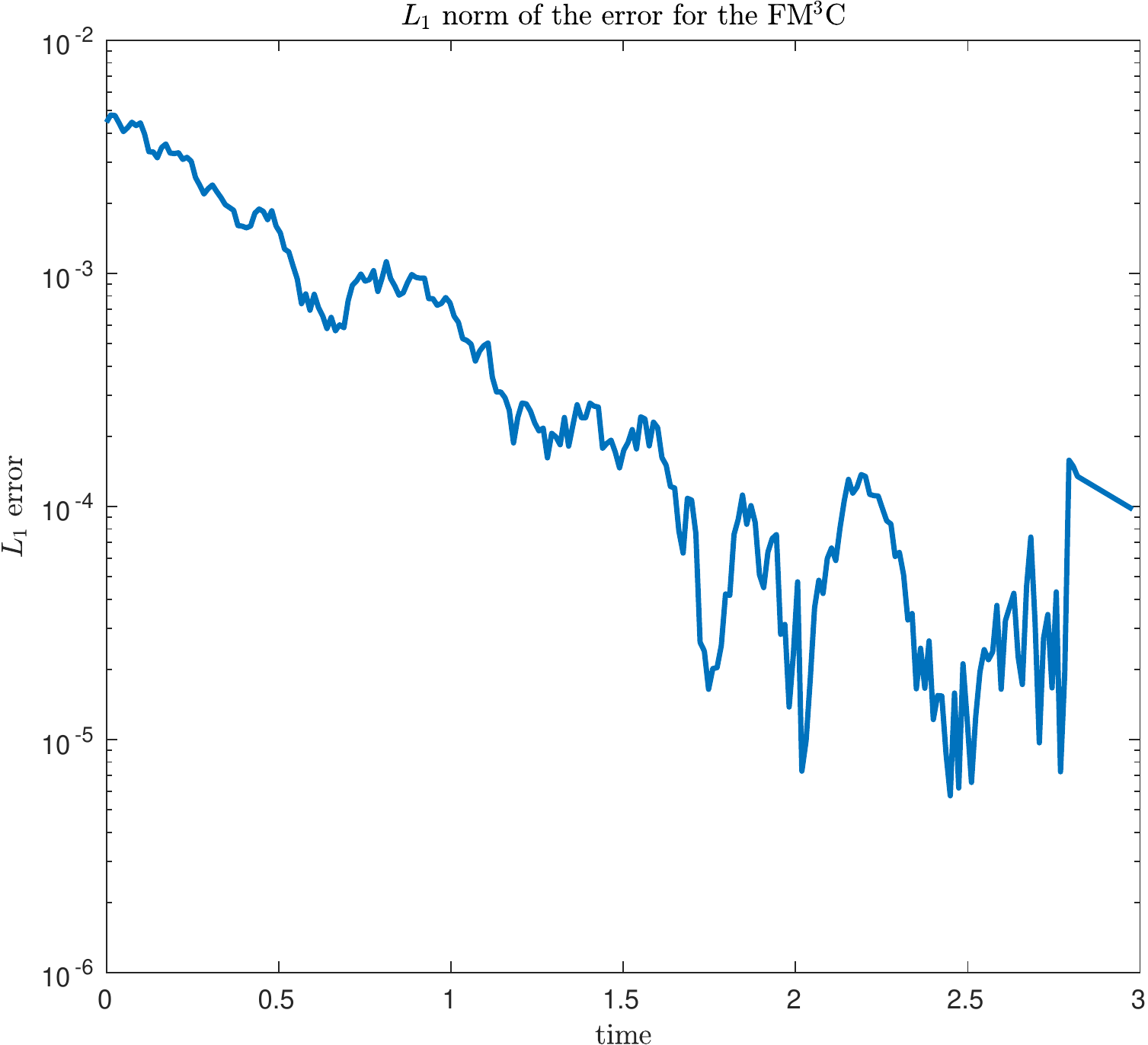}
\includegraphics[scale=0.5]{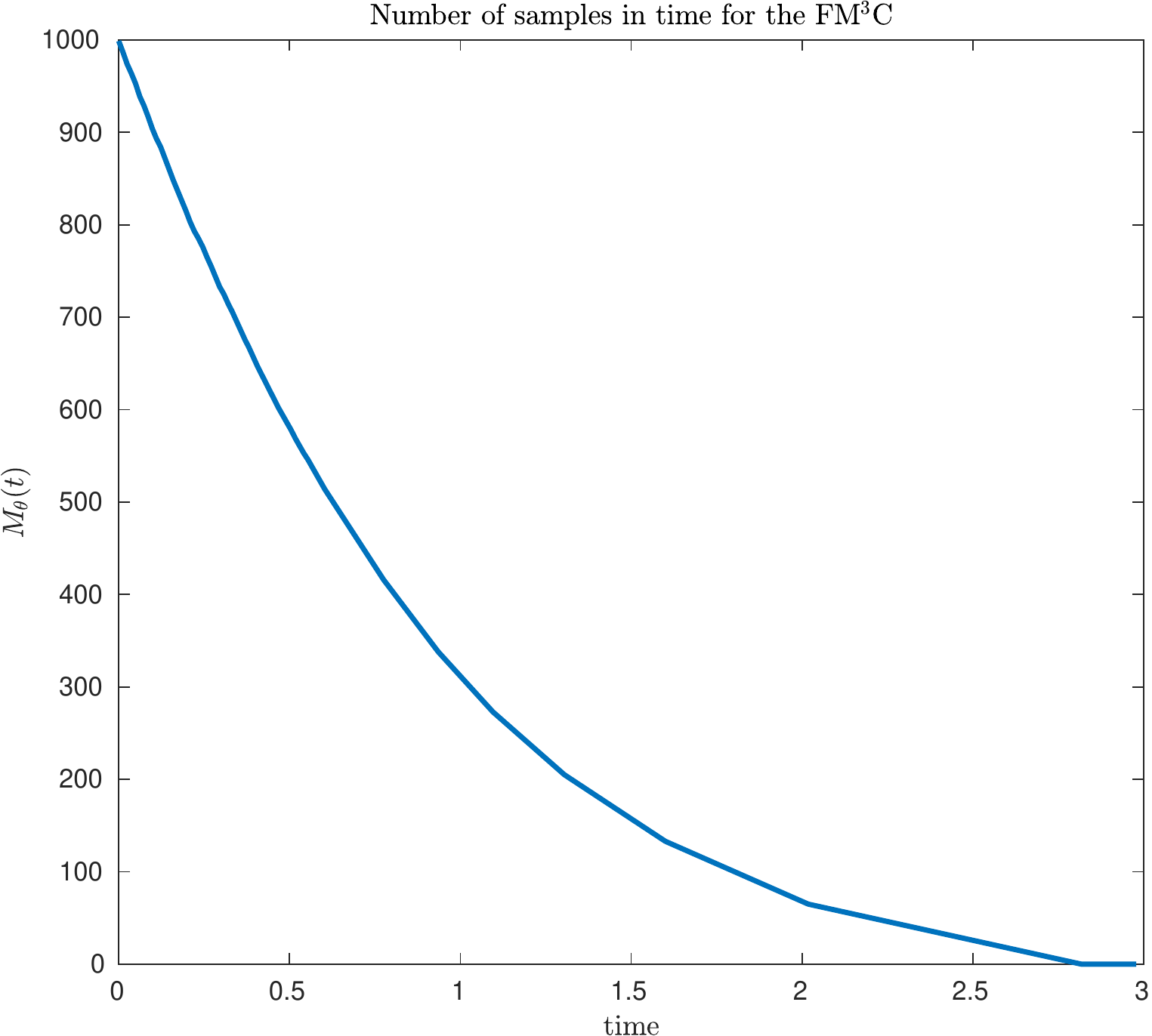}
\caption{Fast M$^3$C method. Left: evolution of the $L_1$ norm of the error for the expected distribution over time. Different lines represent a different number of stochastic inputs over which the expected value is computed. Right: number of random nodes over time. Grid on $[-1,1]$ with $N=100$ points, final time $T=3$ and $\Delta t=\Delta w^2/2$. }
\label{fastMM}
\end{figure}

Finally, in Figure \ref{fastMM} we show the behavior of the fast M$^3$C method. The number of samples for which the time evolution of the perturbation $g$ is considered is reported on the right and it diminishes with time. The corresponding $L_1$ norm of the error is shown on the left. For this case, we increased the initial number of random nodes to $1000$ to highlight the behavior of the fast approach. The number of evolutions of the distribution function computed diminishes exponentially.

\section{Stochastic Galerkin methods}
Among the various methods for UQ in PDEs, stochastic Galerkin (SG) methods based on generalized polynomial chaos (gPC) expansions are very attractive thanks to the spectral convergence property with respect to the random input \cite{CVK,HJX,JXZ,LeMK,PIN_book,PDL,X,XK}. On the other hand, their intrusive nature forces a complete reformulation of the problem and standard schemes for the corresponding deterministic problem cannot be used in a straightforward way.  

In particular, it is well known that, this intrusive formulation may lead to the loss of important structural properties of the original problem, like hyperbolicity, positivity and preservation of large time behavior \cite{CJK,DPL,PIN}. 

In this section we analyze gPC-SG methods for the numerical approximation of stochastic Vlasov-Fokker-Planck equations in the form (\ref{eq:MF_general}). In particular, using the Micro--Macro approach in the gPC-SG setting we show how it is possible to construct methods which preserve the asymptotic behavior of the solution \cite{PZ3}. We mention here related approaches for kinetic equations developed in \cite{HJ,ZJ}.

We recall first some basic notions on Galerkin approximation techniques for stochastic computations.

\subsection{Preliminaries on gPC-SG techniques}\label{sec:gPC_prelim}

Let us consider the function $f(\theta,x,w,t)$, $f\in L^2$ in the random variable $\theta\in\Omega\subseteq \RR$, solution of the differential problem 
\be\label{eq:differential_problem}
\partial_t f(\theta,x,w,t) = \mathcal J[f](\theta,x,w,t),
\ee
with $\mathcal J$ a given differential operator. The present setup of the problem may be naturally extended to a $r-$dimensional vector of random variables $\mathbf{\theta}=(\theta_1,\dots,\theta_r)$. 

We consider the space $\mathbb{P}_M$ of polynomials of degree up to $M$, generated by a family of orthogonal polynomials with respect to the probability density function $p(\theta)$ of the random variable $\theta$, namely $\{\Phi_h(\theta)\}_{h=0}^M$. They form an orthogonal basis of $L^2(\Omega)$, i.e.
\begin{equation}
\mathbb E\Big[\Phi_h(\theta)\Phi_k(\theta)\Big]=\int_{\Omega}\Phi_h(\theta)\Phi_k(\theta)p(\theta)\,d\theta= \mathbb E\Big[\Phi_h^2(\theta)\Big]\delta_{hk}
\end{equation}
where $\delta_{hk}$ is the Kronecker delta function. Let us assume that $p(\theta)$ has finite second order moment, we can represent the function $f(x,w,\theta,t)$ through the complete polynomial chaos expansion as follows
\begin{equation}
f(\theta,x,w,t)=\sum_{m\in\mathbb{N}}\hat{f}_m(x,w,t)\Phi_m(\theta),
\end{equation}
where $\hat{f}_m(x,t)$ is given by 
\begin{equation}\label{eq:hat_fm}
\hat{f}_m(x,w,t)=\mathbb E\Big[f(\theta,x,w,t)\Phi_m(\theta)\Big],\qquad m\in\mathbb{N}.
\end{equation}
The generalized polynomial chaos expansion approximates the solution $f(\theta,x,w,t)$ of \eqref{eq:differential_problem} with its $M$-th order truncation $f^M(\theta,x,w,t)$ and considers the Galerkin projections of differential problem for each $h=0,\dots,M$
\begin{equation}
\partial_t \mathbb E[f(\theta,x,w,t)\cdot \Phi_h(\theta)] =\mathbb E[\mathcal{J}[f](\theta,x,w,t)\cdot \Phi_h(\theta)].
\end{equation} 
Thanks to the orthogonality of the polynomial basis of the space $\mathbb{P}_M$ we obtain a coupled system of $M+1$ purely deterministic equations  
\begin{equation}
\partial_t \hat f_h(x,w,t) = \mathcal J[(\hat f_k)_{k=0}^M](x,w,t),\qquad h=0,\dots,M.
\end{equation}
These subproblems must then be solved through suitable numerical techniques. The approximation of the statistical quantities of interest are defined in terms of the introduced projections. From \eqref{eq:hat_fm} being $\Phi_0\equiv 1$ we have
\begin{equation}\begin{split}
\mathbb E[f(\theta,x,w,t)] &= \hat f_0(x,w,t),
\end{split}\end{equation}
 and thanks to the orthogonality it is possible to show that 
 \begin{equation}
 \begin{split}
 Var[f(\theta,x,w,t)] &= \mathbb E\Bigg[ \Big(\sum_{h=0}^M \hat f_h \Phi_h(\theta)-\hat f_0\Big)^2 \Bigg] \\
 &=\sum_{h=0}^M \hat f_h^2(x,w,t)\mathbb E[\Phi_h^2(\theta)]-\hat f_0^2(x,w,t). 
\end{split}\end{equation}

\subsubsection{gPC-SG methods for Vlasov--Fokker--Planck equations}\label{sec:gPC_MF}
Let us consider the stochastic Vlasov--Fokker--Planck equation \eqref{eq:MF_general} with a nonlocal drift $\mathcal B[\cdot]$ of the form (\ref{eq:Bf}).

The gPC-SG approximation is given by the following system of deterministic differential equations
\be\begin{split}\label{eq:system_gPC}
\partial_t \hat f_h(x,w,t)+&\mathcal L[\hat f_h(x,w,t)]=\\
&\nabla_w\cdot \left[ \sum_{k=0}^M b_{hk}[\hat f](x,w,t)\hat f_k(x,w,t)+\nabla_w D(x,w)\hat f_h(x,w,t) \right],
\end{split}\ee
where  
\be\label{eq:system_gPC2}
b_{hk}[\hat f](x,w,t) = \dfrac{1}{\|\Phi_h \|^2_{L^2}} \sum_{m=0}^M\int_{\Omega} \mathcal B[\hat f_m]\Phi_k(\theta)\Phi_m(\theta) dp(\theta).
\ee
Note that, due to the nonlinearity of Fokker--Planck problems, we obtain a coupled system of deterministic Vlasov--Fokker--Planck equations describing the evolution of each projection. In vector notations we have 
\be
\partial_t \hat{\mathbf{f}}(x,w,t) +\mathcal L [\hat{\mathbf{f}}](x,w,t)= \nabla_w \cdot [\mathbf{B}[\hat{\mathbf f}](x,w,t)\hat{\mathbf{f}}(x,w,t)+\nabla_w \mathbf{D}(x,w)\hat{\mathbf{f}}(x,w,t) ],
\ee
where $\hat{\mathbf{f}}=(\hat f_0,\dots,\hat f_M)^T$ and the component of the $(M+1)\times (M+1)$ matrix $\mathbf{B}[\hat{\mathbf{f}}](x,w,t)$ are given by \eqref{eq:system_gPC2}.

In a similar way, we can derive the gPC-SG formulation of stochastic Vlasov--Fokker--Planck equations with uncertain diffusion terms.
 
\begin{remark}
In case the uncertainty is present only in the initial data, and therefore $\mathcal B[f](x,w,\theta,t)=\mathcal B(x,w,t)$, the matrix $\mathbf B$ is diagonal and we need to solve the decoupled system of Vlasov type equations
\be
\partial_t \hat f_h(x,w,t)+ \mathcal L[\hat f_h](x,w,t) = \nabla_w \cdot [b_{hh}\hat f_h(x,w,t) + \nabla_w D(x,w) \hat f_h(x,w,t)],
\ee
$h=0,\dots,M$. Hence, a structure preserving approach as in Section \ref{sec:structure_preserving} may be introduced in order to preserve the large time behavior of the collision step of each projection by defining a family of weight functions 
\be\begin{split}
\lambda_{i+1/2}^h &= \dfrac{D_{i+1/2}}{\Delta w}\int_{w_i}^{w_{i+1}}\dfrac{b_{hh}(x,w,t)+D'(x,w)}{D(x,w)}dw, \\
\delta_{i+1/2}^h &= \dfrac{1}{\lambda_{i+1/2}^h}+\dfrac{1}{1-\exp(\lambda_{i+1/2}^h)}.
\end{split}\ee
In this setting the scheme capture with arbitrary accuracy the steady state and the expected value of the numerical solution is kept nonnegative. However, for more general nonlocal type operators $\mathcal B[\cdot]$ this approach cannot be applied for the construction of a stochastic Galerkin expansion which preserves the steady state solution and nonnegativity of the mean. 
\end{remark}

\subsection{A Micro--Macro gPC approach}\label{sec:MM_gPC}
We discussed in the previous section how the gPC-SG method for stochastic Fokker--Planck equations generates a coupled system of partial differential equations. Although gPC-SG guarantees spectral convergence on the random field under suitable regularity assumptions, its accuracy in  describing the long--time solutions of the problems is limited and depends on the particular scheme for solving the coupled system.

Let us consider suitable regularity assumptions on the initial distribution such that the stochastic Fokker--Planck problem admits the unique steady state solution $f^{\infty}(\theta,w)$. With the aim of preserving the steady states of the problem in the Galerkin setting we  introduce a Micro--Macro gPC-SG scheme. Thanks to the formalism introduced in Section \ref{sec:gPC_prelim} and by analogy with \eqref{eq:f_MM} the Micro--Macro gPC decomposition for all $M\ge 0$ reads \cite{PZ3}
\be\label{eq:finf_g_gPCMM}
f^M(\theta,w,t) = f^{\infty,M}(\theta,w)+g^M(\theta,w,t),\qquad w\in\RR^{d_w}, t\ge 0,
\ee
where 
\[
f^{\infty,M}(w,\theta) = \sum_{h=0}^M \widehat{f^{\infty}}_h(w)\Phi_h(\theta), \qquad \widehat{f^{\infty}}_h(w) = \int_{\Omega}f^{\infty}(\theta,w)\Phi_h(\theta)dp(\theta).
\]
Being equation \eqref{eq:finf_g_gPCMM} equivalent to require $\hat{\mathbf{f}}=\widehat{\mathbf{f}^{\infty}}+\hat{\mathbf{g}}$, we can reformulate the original problem in terms of $\hat{\mathbf{g}}$. Equation \eqref{eq:system_gPC} may be reformulated for all $h=0,\dots,M$ in terms of the nonequilibrium part of the Micro--Macro gPC decomposition $\hat g_h$ as follows
\be
\label{eq:sfp}
\begin{cases}
\partial_t \hat{g}_h(w,t) &= \hat{\mathcal J}_h(\hat g,\hat g)(w,t)+ \hat{\mathcal N}_h(\widehat{f^{\infty}},\hat g)(w,t),\\
f^M(w,\theta,t) &= f^{\infty,M}(w,\theta)+g^M(w,\theta,t),
\end{cases}
\ee
where the operator $\hat{\mathcal J}_h$ is the Galerkin projection of the quadratic operators of the collisional type defined in \eqref{eq:J_def} and $\hat{\mathcal N}_h$ is a linear operator defined as 
\be\begin{split}\label{eq:JN_gPCMM}
\hat{\mathcal J}_h(\hat g,\hat g)(w,t) &= \nabla_w\cdot \Big[ \sum_{k=0}^M b_{hk}[\hat g]\hat g_k(w,t)+\nabla_w D(w)\hat g_h(w,t) \Big],\\
\hat{\mathcal N}_h(\widehat{f^{\infty}},\hat g)(w,t) &= \nabla_w \cdot \Big[ \sum_{k=0}^M b_{hk}[\widehat{ f^{\infty}}]\hat g_k(w,t) + b_{hk}[\hat g]\widehat{f^{\infty}}_k(w) \Big]. 
\end{split}\ee

Now, the equilibrium state of each gPC projection is $\hat g_h\equiv 0$ and any consistent schemes for the numerical approximations of the differential terms in (\ref{eq:JN_gPCMM}) admits $\hat g_h \equiv 0$ as equilibrium state for all $h=0,\dots,M$. For example, we can use a standard central difference approximation scheme for the differential terms in (\ref{eq:JN_gPCMM}) to achieve second order accuracy for transient times and exact preservation of the steady state asymptotically.

\subsection{Numerical results}
\begin{figure}[t]
\centering
\includegraphics[scale=0.5]{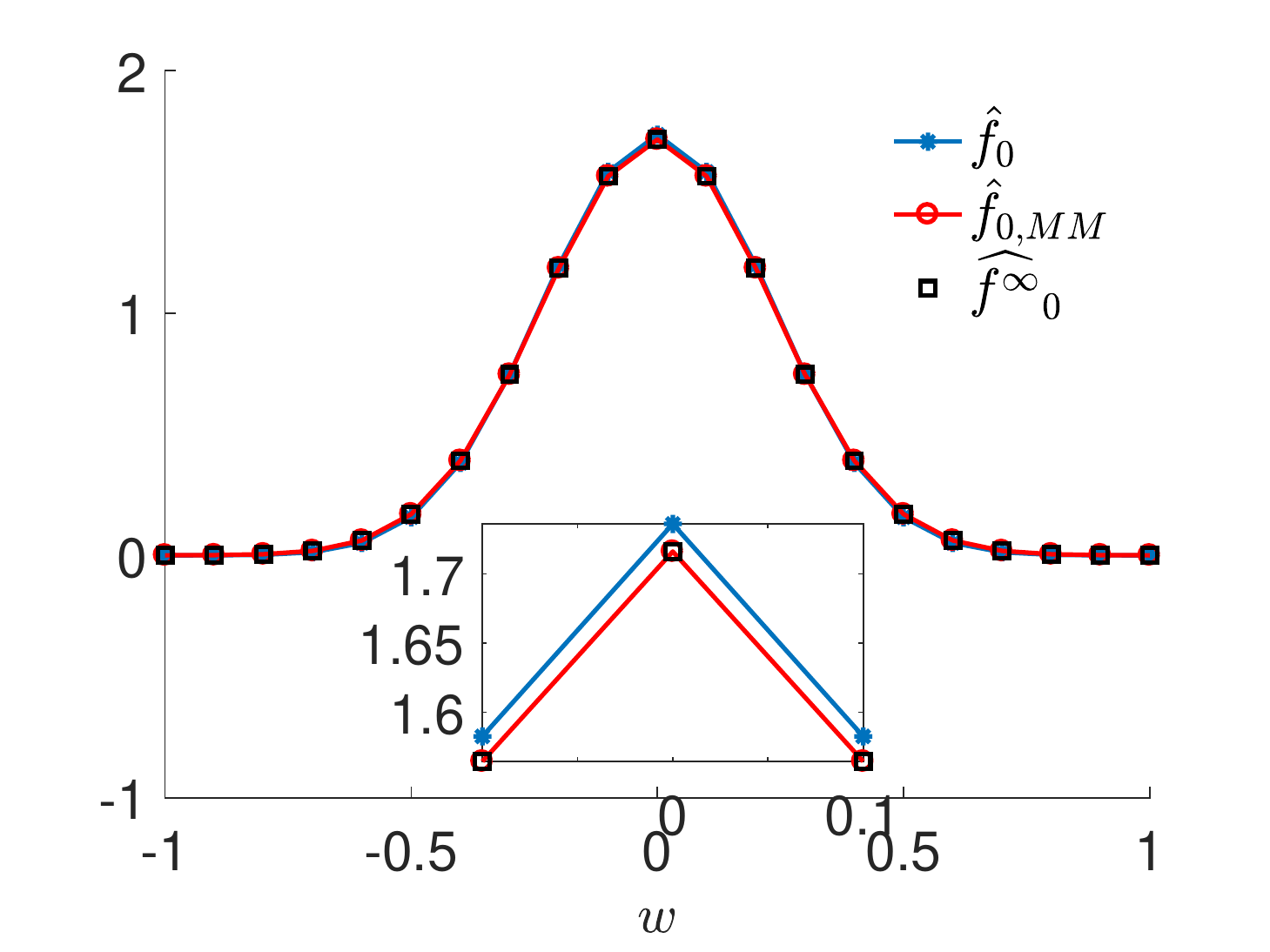}
\includegraphics[scale=0.5]{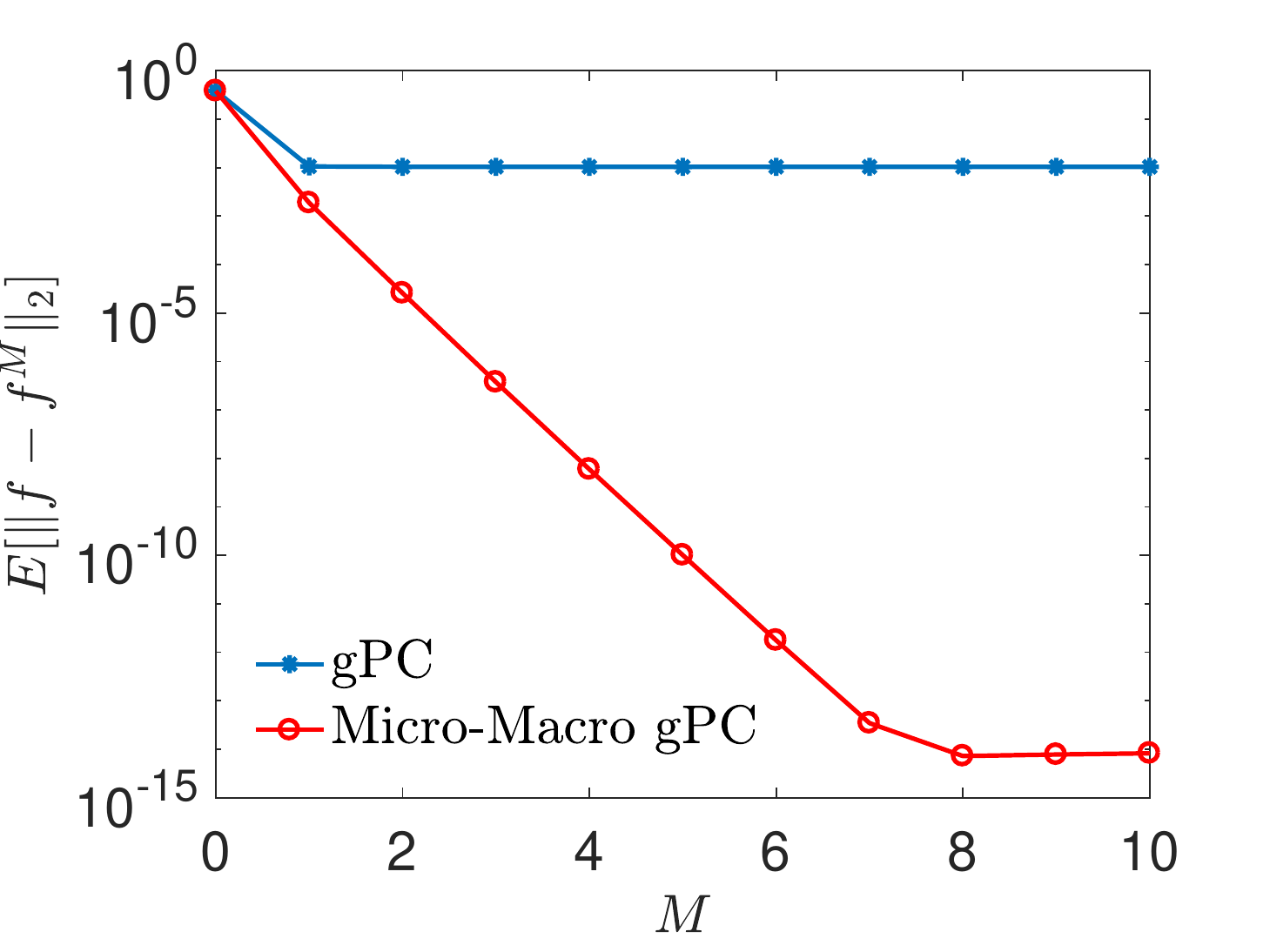}
\caption{Left: expected solution of the uncertain Fokker--Planck equation obtained through standard SG and MM-SG methods with central differences and $M=10$. Right: expected $L^2$ error for standard SG and MM-SG. In both cases we considered a discretization of the interval $[-1,1]$ with $N=21$ gridpoints and $\Delta t=\Delta w/2$. }
\label{fig:gPC_FP}
\end{figure}

We consider the evolution of the Fokker--Planck equation \eqref{eq:FP_num_SP} with the uncertain initial condition \eqref{eq:FP_initial}. 
Following the set--up introduced in the previous section, we obtain the SG system of equations
\be\label{eq:F_gPC_system}
\partial_t \hat f_h(w,t) = \partial_w \Big[ w \hat f_h(w,t) + \partial_w \sum_{k = 0}^M d_{hk}\hat f_k(w,t) \Big],
\ee
with 
\[
d_{hk} = \dfrac{1}{\| \Phi_h \|^2_2} \int_{\Omega}  {T(\theta)}\Phi_h(\theta)\Phi_k(\theta)dp(\theta).
\]
In order to build the Micro--Macro gPC decomposition of the SG system we take advantage of the analytical solution given by the Maxwellian distribution \eqref{eq:maxw}, which can be approximated by its $M$--order truncation as in Section \ref{sec:MM_gPC}. Therefore we aim at solving the modified problem for all $h=0,\dots,M$
\begin{equation}
\begin{cases}
\partial_t \hat g_h(w,t) = \partial_w \Big[w\hat g_h(w,t)+\partial_w \sum_{k=0}^M d_{hk}\hat g_k(w,t)  \Big], \\
f^M(\theta,w,t)=g^M(\theta,w,t)+f^{\infty,M}(\theta,w).
\end{cases}
\end{equation}
In all our numerical examples we use second order central difference approximations of the derivatives in $w$.
In Figure \ref{fig:gPC_FP} we compare the numerical long time solution obtained through a standard SG system  \eqref{eq:F_gPC_system} and the Micro--Macro SG system (MM). We can observe how the Micro--Macro gPC-SG method gives an accurate description of the expected steady state of the problem, on the contrary the error of the standard gPC-SG method saturates at the accuracy obtained with the central differences. 

\section{Other applications}
\label{sec:6}
In this section we present several numerical examples of stochastic Fokker--Planck and Vlasov--Fokker-Planck equations solved with the schemes introduced in the previous sections. In particular we focus on some recent models in socio--economic and life sciences as discussed in the Introduction.

\begin{figure}[t]
\centering
\includegraphics[scale=0.5]{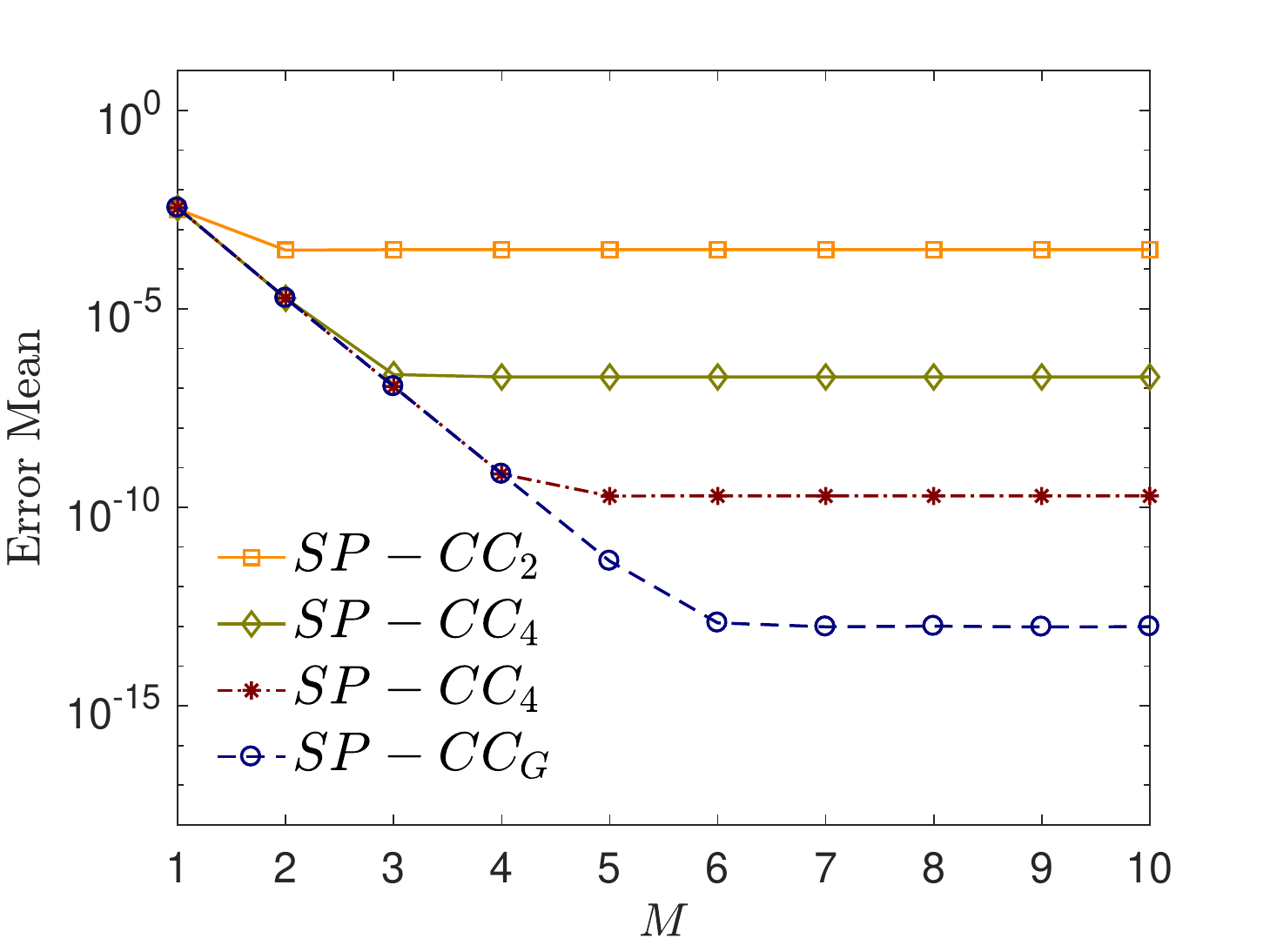}
\includegraphics[scale=0.5]{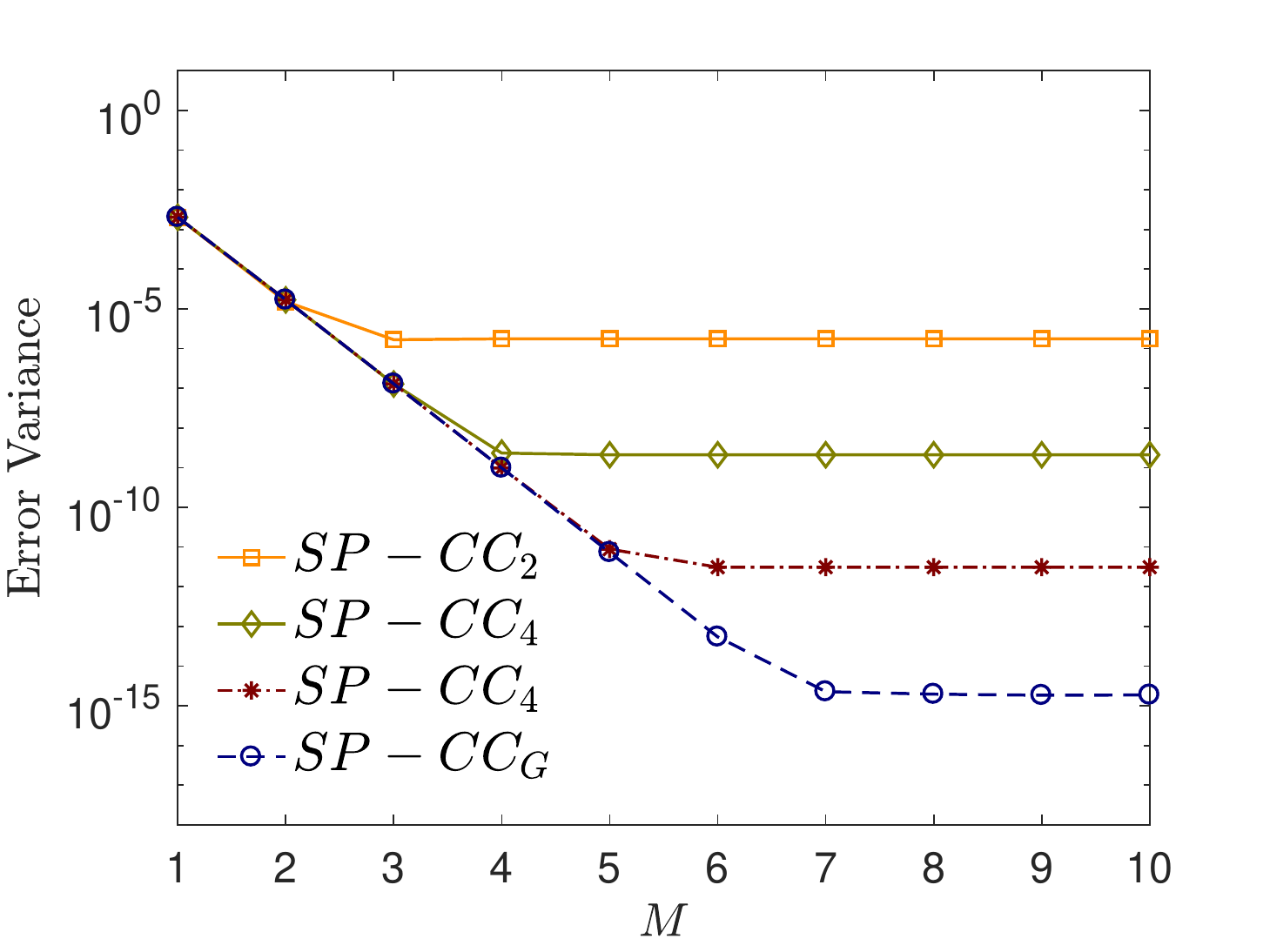}
\caption{Example 1. Values of $L^1$ error in the estimation of the expected solution (left) and its variance (right) for  $T=20$ and for an increasing number of collocation nodes. The numerical error has been computed with respect to the expected analytical solution (left) and its variance (right), see \eqref{eq:exact_ss_op}. We compare the error for the SP--CC scheme with different quadrature methods in case of random interaction $P(\theta)=0.75+\theta/4$, $\theta\sim U([-1,1])$. Initial distribution \eqref{eq:initial_data_opinion}, $\sigma^2/2=0.1$, $N = 80$, $\Delta t=\Delta w^2/(2\sigma^2)$.  }
\label{fig:collocation_error_op}
\end{figure}

\subsection{Example 1: Opinion model with uncertain interactions}
Let us consider a distribution function $f=f(\theta,w,t)$ describing the density of agents with opinion $w\in I=[-1,1]$ whose evolution is given in terms of a stochastic Fokker-Planck equation characterized by the nonlocal term \eqref{eq:BD_opinion} with uncertain compromise propensity function $P(\theta,w,_*)\in [0,1]$. In the following we will solve the problem both in the collocation and in the Galerkin setting. 

We consider as deterministic initial distribution $f(\theta_k,w,0)=f_0(w)$ for all $k=1,\dots,M$, with
\be\label{eq:initial_data_opinion}
f_0(w,0) = \beta \left[ \exp(-c(w+1/2)^2)+\exp(-c(w-1/2)^2) \right],\qquad c=30,
\ee
with $\beta>0$ a normalization constant and let $u = \int_{-1}^1 wf_0(w)dw$ the mean opinion. We choose a uniformly distributed random input $\theta\sim U([-1,1])$ and a random interaction function of the form $P(\theta)=0.75+\theta/4$. 

We discretize the random variable by considering the first $M>1$ Gauss--Legendre collocation nodes.  In Figure \ref{fig:collocation_error_op} we compute the relative $L_1$ error for mean and variance with respect to the exact steady state \eqref{eq:exact_ss_op} using $N=80$ points for the $SP-CC$ scheme with various quadrature rules adopted for the evaluation of the weights function in (\ref{eq:lambda_high}). Singularities at the boundaries in the integration of \eqref{eq:lambda_high} can be avoided using open Newton--Cotes methods. In the sequel, we will adopt the notation $SP-CC_k$, $k=2,4,6,G$, to denote the structure preserving schemes with Chang--Cooper flux when \eqref{eq:lambda_high} is approximated with second, fourth, sixth order open Newton--Cotes or Gaussian quadrature, respectively. 

\begin{figure}
\centering
\includegraphics[scale=0.4]{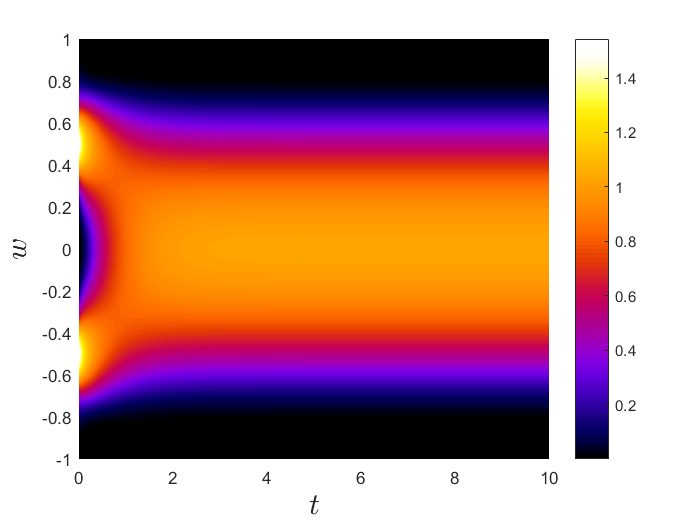}
\includegraphics[scale=0.4]{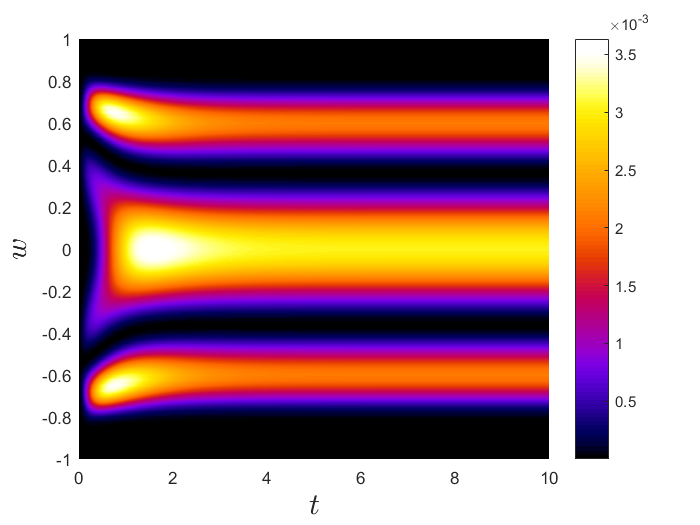}
\caption{Example 1. Evolution of $\mathbb E^M[f(w,\theta,t)]$ (left) and $Var^M[f(w,\theta,t)]$ (right) for the opinion model obtained with $M=10$ collocation points and the $SP-CC_G$ scheme over the time interval $[0,10]$. }
\label{fig:opinion_evo_meanvar}
\end{figure}

In Figure \ref{fig:opinion_evo_meanvar} the time evolution of the expected solution and variance are given. We can observe from the estimation of the variance the regions of higher variability of the  expected solution due to uncertain interactions.  The evolutions of the statistical quantities have been computed through a collocation $SP-CC_G$ method with $20$ quadrature points for the evaluation of \eqref{eq:lambda_high}. 

\begin{figure}
\centering
\includegraphics[scale=0.5]{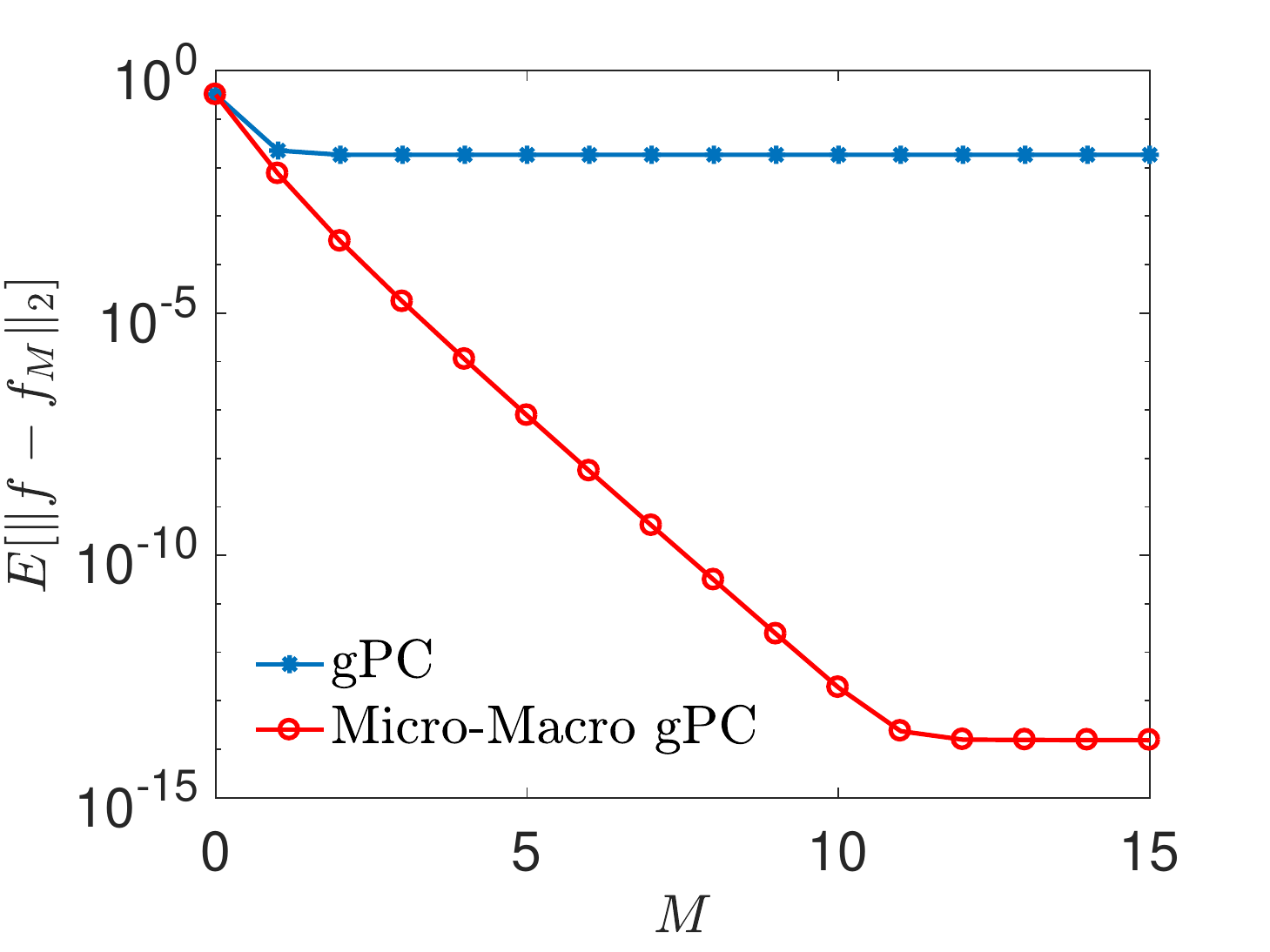}
\includegraphics[scale=0.5]{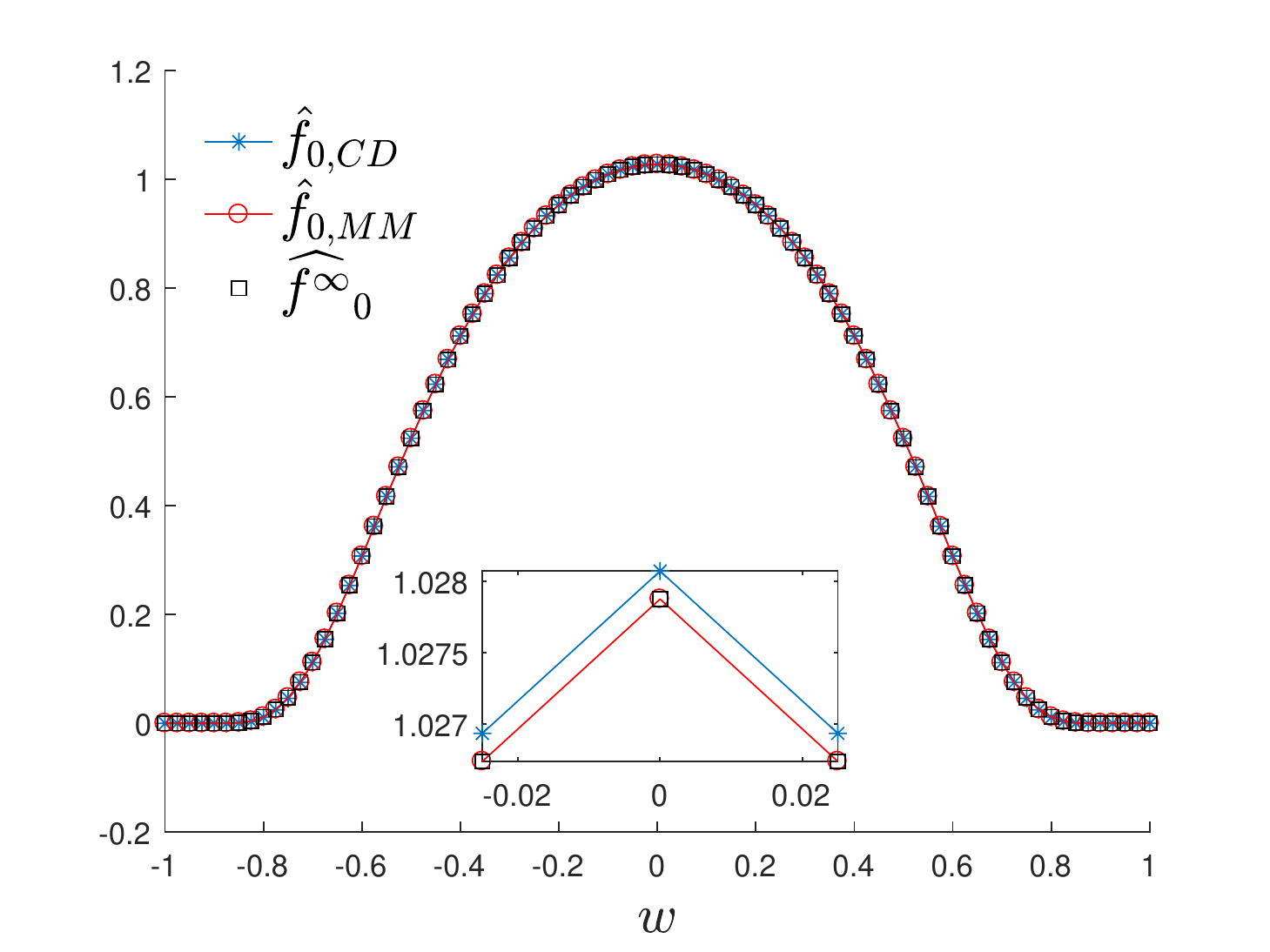}
\caption{Example 1. Left: Estimation of $\mathbb E[\|f-f^M \|_{L^2}]$ computed at time $t=20$ and for an increasing number of $M\ge 0$, we compare the errors computed through a standard method for the solution of the system of coupled PDEs of the Fokker--Planck type, with the Micro--Macro gPC method. We used $N = 80$ gridpoints, $\sigma^2/2=0.1$, $\Delta t^2 = \Delta w^2/(2\sigma^2)$.   Right: Large time behavior for the estimated expected solution of the opinion model. We can observe how the Micro--Macro gPC method is able to capture with high accuracy the steady expected solution of the problem.   }
\label{fig:opinion_gpc}
\end{figure}

Finally, as in Section \ref{sec:MM_gPC} we consider a Micro-Macro gPC Galerkin setting based on the knowledge of the stationary solution (\ref{eq:exact_ss_op}).  

In Figure \ref{fig:opinion_gpc} we present the behavior numerical error $\mathbb E[\|f^{\infty}-f^M \|_2]$ for large time where the differential terms in $w$ are solved by central differences. We report also the large time behavior for the expected solution in both schemes, where it is possible to observe how the Micro--Macro gPC is able to capture with high accuracy the steady state of the problem. 

\begin{figure}
\centering
\includegraphics[scale=0.5]{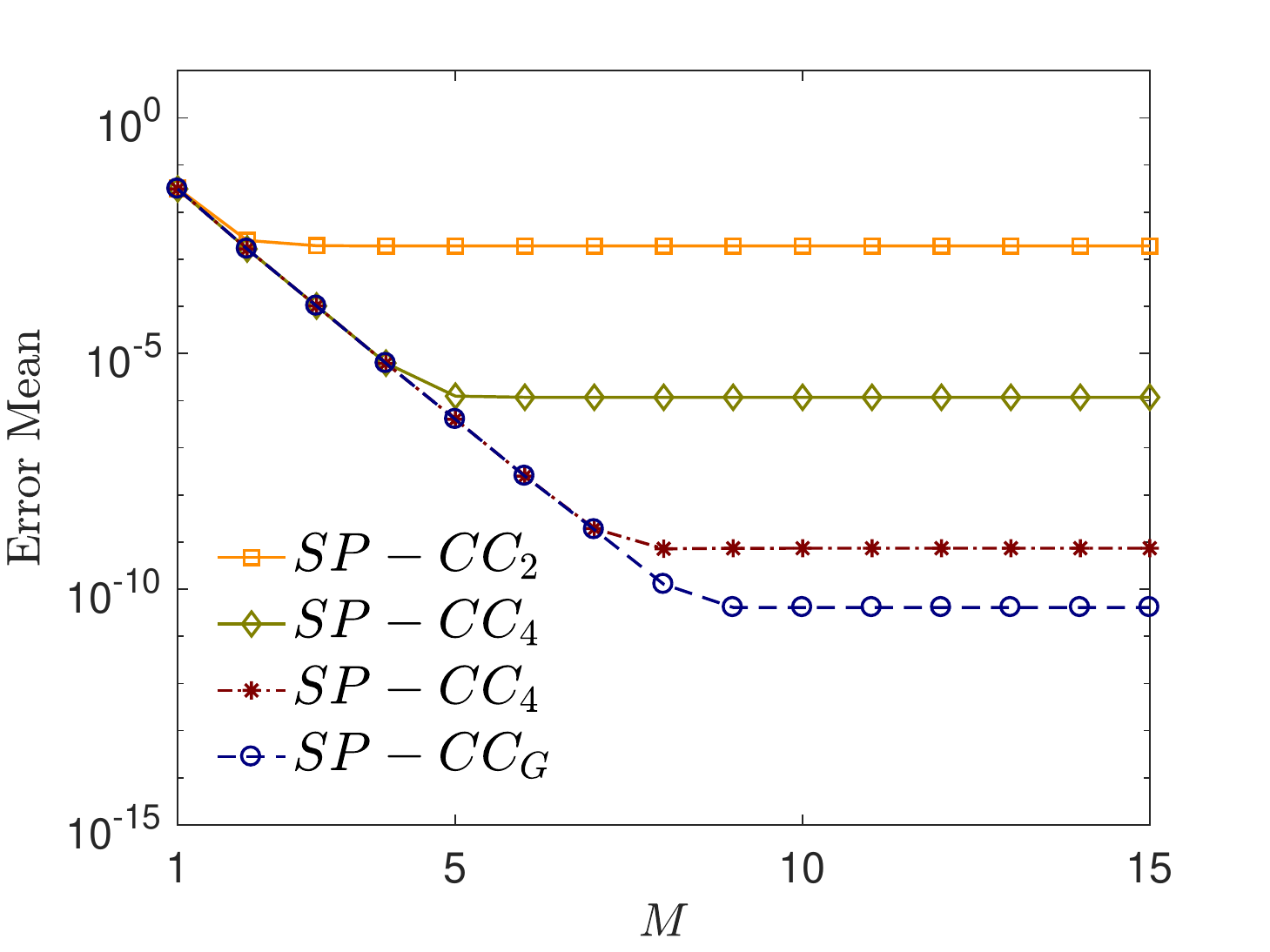}
\includegraphics[scale=0.5]{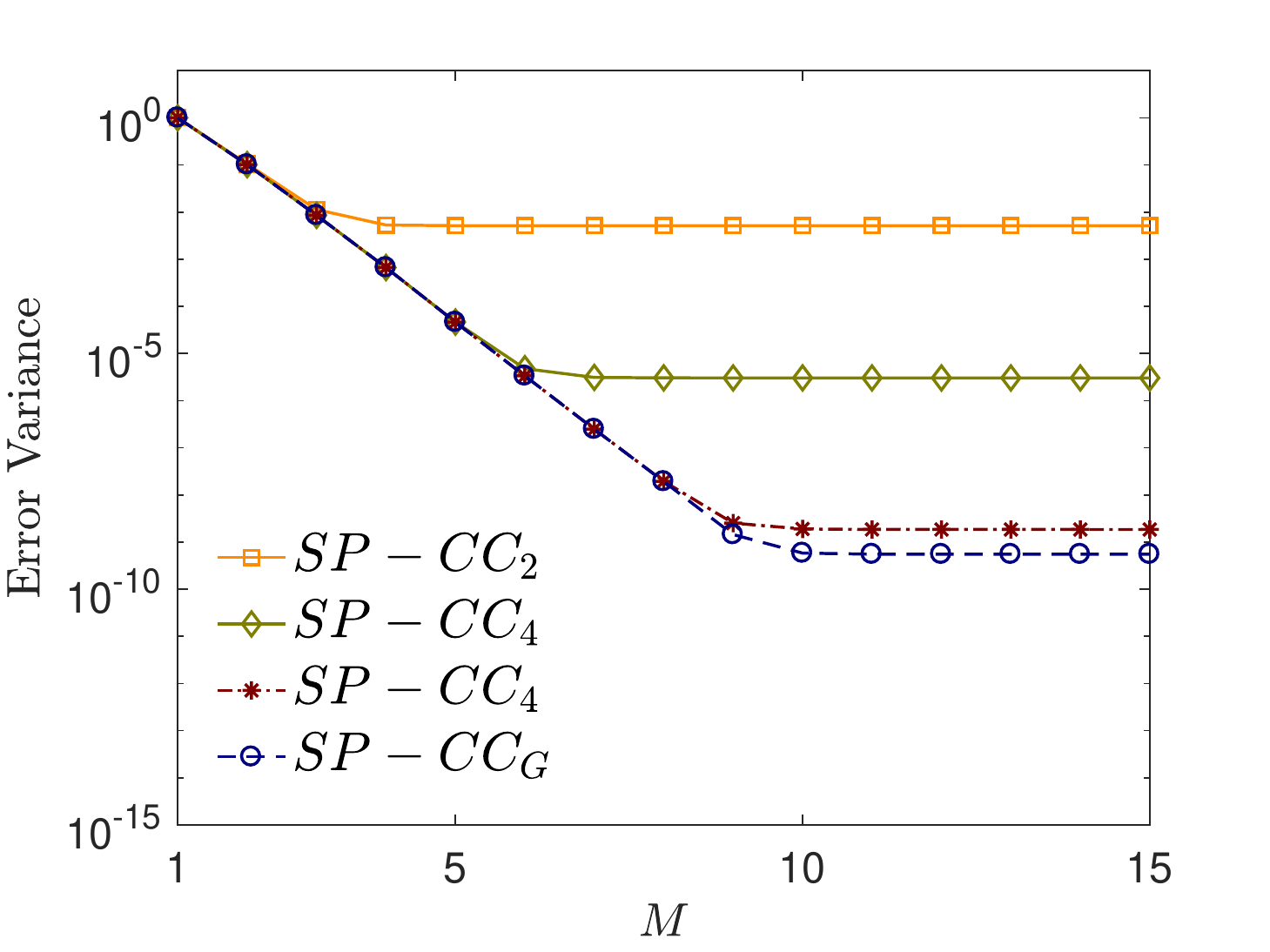}
\caption{Example 2. Error for the SP--CC scheme with different quadrature methods in case of random diffusion constant $\sigma^2(\theta)=0.1+\theta/200$, $\theta\sim U([-1,1])$. We report the $L^1$ relative error in the estimation of the expected solution (left) and its variance (right) for $T=20$ and for an increasing number of nodes in the random space. The  numerical error has been computed with respect to the expected analytical solution (left) and its variance (right) obtained from \eqref{eq:exact_ss_we}. The initial distribution $f_0(w)$ is \eqref{eq:initial_data_wealth}, we consider the domain $[0,L]$, $L=10$ with $N = 200$ points and $\Delta t= \Delta w/L$ with a semi-implicit approximation.}
\label{fig:collocation_error_finance}
\end{figure}

\subsection{Example 2: Wealth evolution with uncertain diffusion}
We consider the Fokker-Planck equation defined by \eqref{eq:BD_wealth} where now $f=f(\theta,w,t)$ with $w\in\RR^+$ represent the wealth of the agents and the uncertainty acts on the diffusion parameter. We consider the deterministic initial distribution $f(\theta_k,w,0)=f_0(w)$ for all $k=1,\dots,M$ with
\be\label{eq:initial_data_wealth}
f_0(w,0) = \beta \exp \Big\{ -c(w-\tilde u)^2 \Big\}, \qquad c = 20,\quad \tilde u = 2,
\ee
where $\beta>0$ is a normalization constant. To deal with the truncation of the computational domain in the interval $[0,L]$, following \cite{PZ1}, after introducing $N$ grid points we consider the quasi stationary boundary condition in order to evaluate $f_{N}(\theta_k,t)$, i.e.
\be\begin{split}\label{eq:finance_SS}
\dfrac{f_{N}(\theta_k,t)}{f_{N-1}(\theta_k,t)} = \exp\Big\{ -\int_{w_{N-1}}^{w_{N}} \dfrac{\mathcal B[f](\theta_k,w,t)+D' (\theta_k,w)}{D(\theta_k,w)}dw \Big\},
\end{split}\ee
for all $k=1,\dots,M$. In Figure \ref{fig:collocation_error_finance} we report in a semilog scale the relative $L_1$ error for mean and variance with respect to the exact steady state introduced in \eqref{eq:exact_ss_we} of the semi--implicit SP--CC  scheme for several integration methods with $N = 200$ points over with $L=10$, and an increasing number of collocation nodes $M=1,\dots,15$. The time step is chosen in such a way that the CFL condition for the positivity of the semi--implicit scheme is satisfied, i.e. $\Delta t=O(\Delta w)$ see Section \ref{sec:prop_SP}. For the tests we considered $\sigma^2(\theta)=0.1+5\times 10^{-2}\theta$, where $\theta\sim U([-1,1])$. We can observe how the error decays exponentially for an increasing number of collocation nodes.

\begin{figure}
\centering
\includegraphics[scale=0.4]{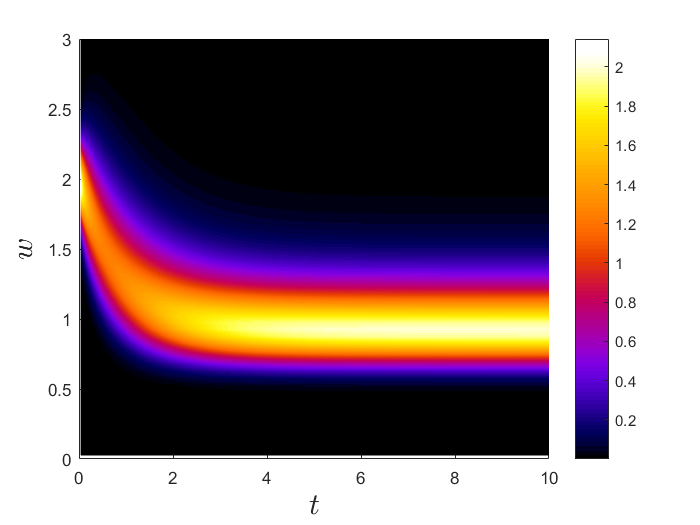}
\includegraphics[scale=0.4]{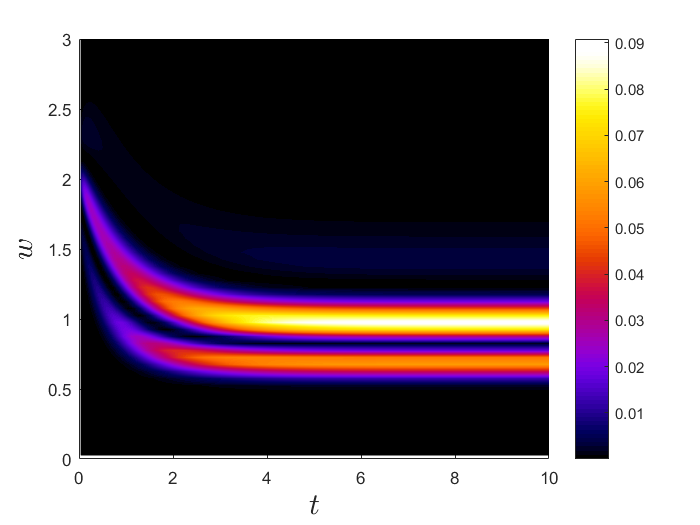}
\caption{Example 2. Evolution of expected solution $\mathbb E^M[f(\theta,w,t)]$ (left) and its variance $\textrm{Var}^M[f(\theta,w,t)]$ (right) for the wealth evolution model. The evolution is computed through $M=10$ collocation points and the $SP-CC_G$ scheme over the time interval $[0,10]$, $\Delta t= \Delta w/L$ with $w\in[0,L]$, $L=10$. }
\label{fig:wealth_evo_meanvar}
\end{figure}

\begin{figure}
\centering
\includegraphics[scale=0.5]{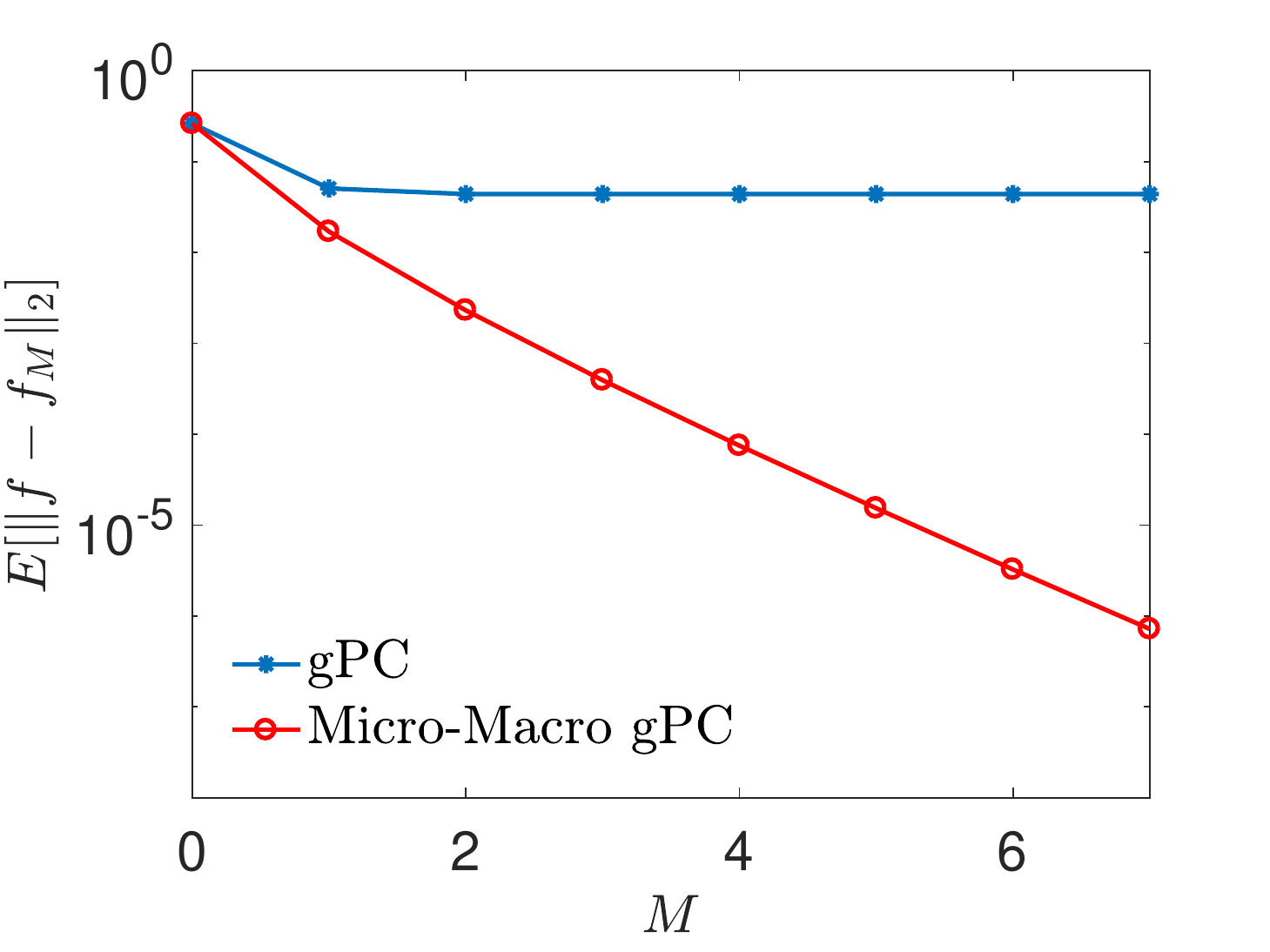}
\includegraphics[scale=0.5]{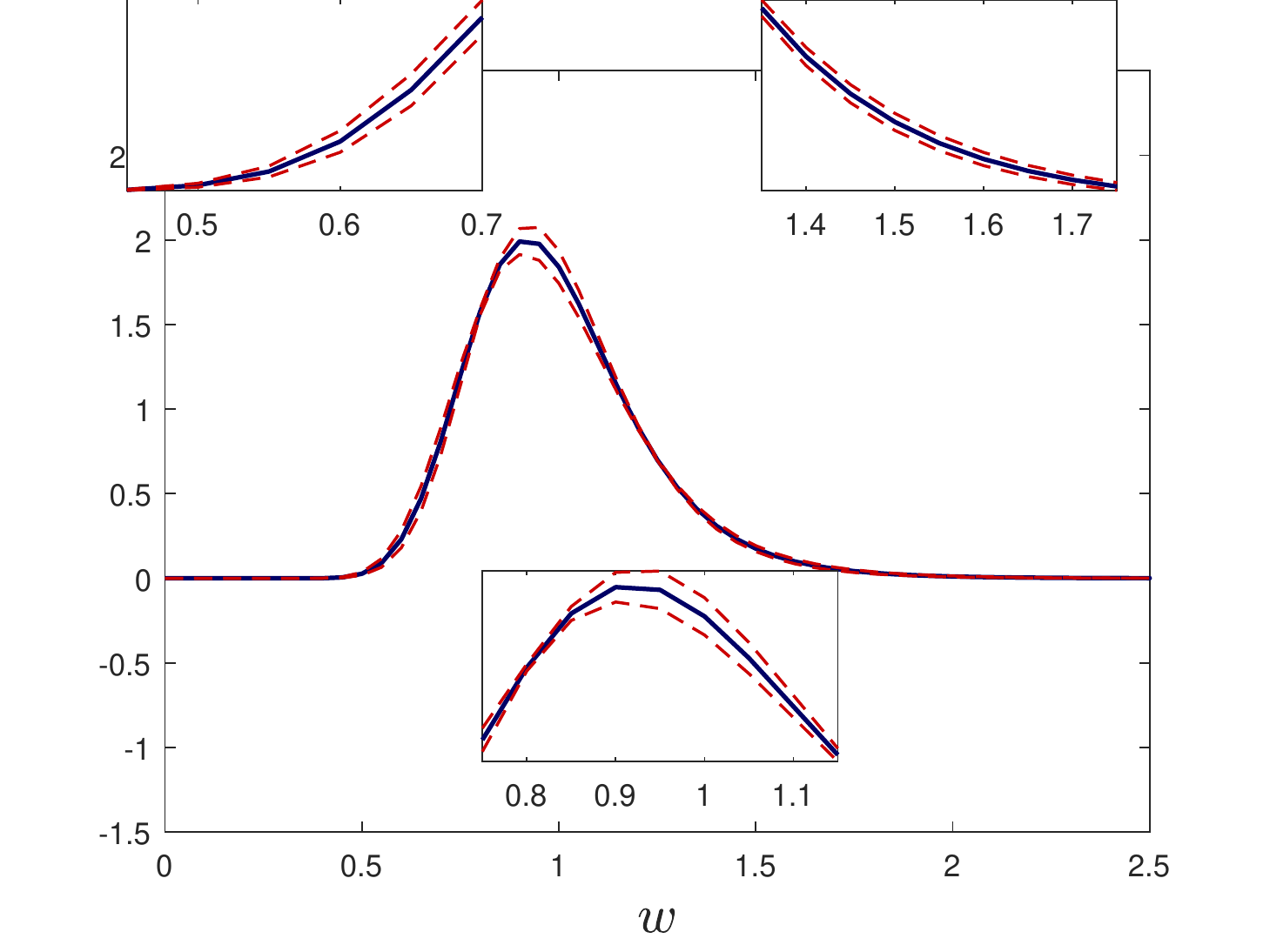}
\caption{Example 2. Left: Estimation of $\mathbb E[\|f-f^M \|_{L^2(\Omega)}]$ computed at time $T=20$ and for an increasing number of $M\ge 0$, we compare the errors computed through a standard gPC-SG method and the Micro--Macro gPC-SG method. We used $N = 200$ gridpoints, $\sigma^2=0.1 + \theta/200$.   Right: Statistical dispersion of the expected asymptotic solution of the wealth distribution model calculated with the Micro--Macro gPC-SG method.}
\label{fig:wealth_MMgpc}
\end{figure}

Next we consider the SG-gPC formulation of the equation for the wealth evolution. Since in this case the uncertainty enters in the definition of the diffusion variable $\sigma^2=\sigma^2(\theta)$ taking $a(\cdot,\cdot)\equiv 1$ the analytical steady state solution of the problem is given in \eqref{eq:exact_ss_we} and we can consider the Micro--Macro gPC scheme as in Section \ref{sec:MM_gPC}.

In Figure \ref{fig:wealth_MMgpc} we compare the error for a standard gPC approximation and the Micro--Macro gPC. In both cases central differences have been used for the differential terms in $w$. We can see how $\mathbb E[\| f-f^M \|_2]$ computed at time $T=20$, close to the stationary solution, decreases in relation to the number of terms of the gPC approximation whereas the standard gPC show a limited accuracy given by the error in approximating the large time behavior of the problem.

\begin{figure}[t]
\centering
\subfigure[$t=1$, MC]{\includegraphics[scale=0.32]{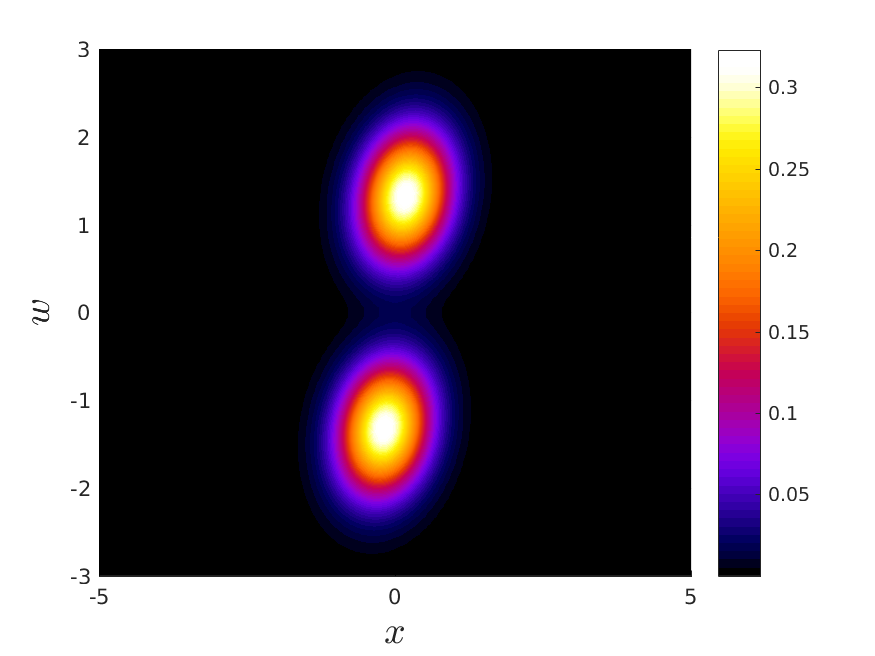}}
\subfigure[$t=3$, MC]{\includegraphics[scale=0.32]{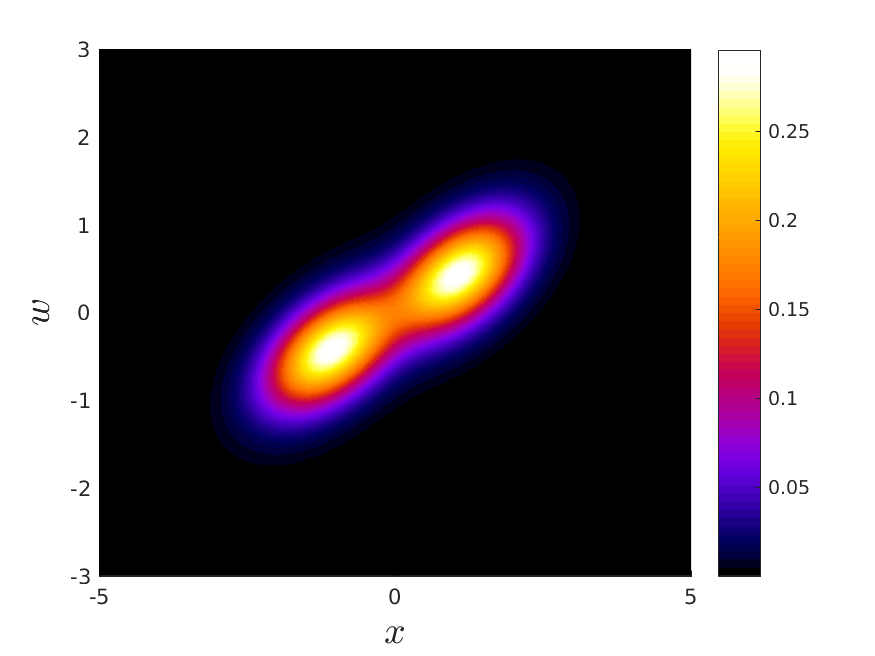}}
\subfigure[$t=6$, MC]{\includegraphics[scale=0.32]{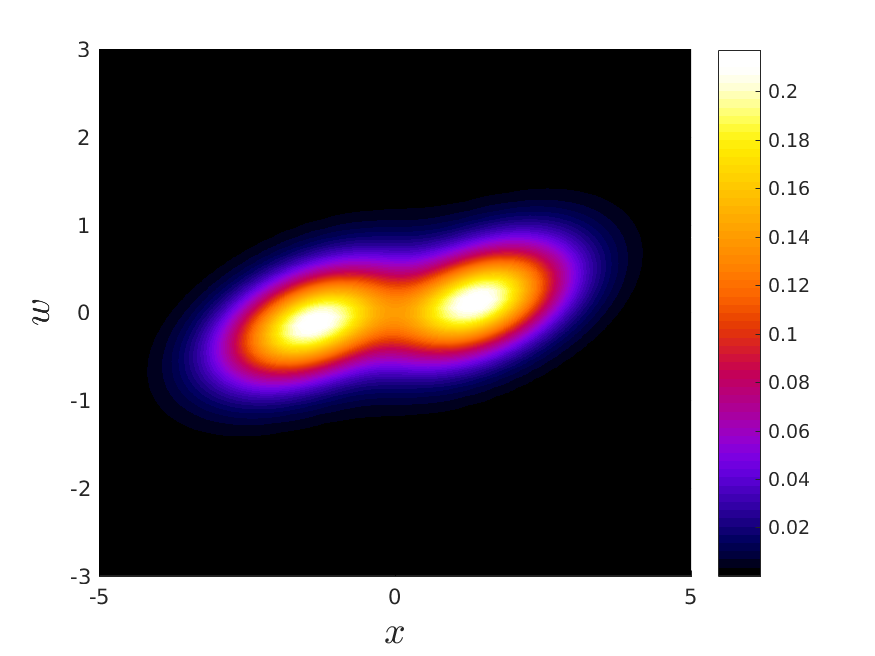}}\\

\subfigure[$t=1$, M$^3$C]{\includegraphics[scale=0.32]{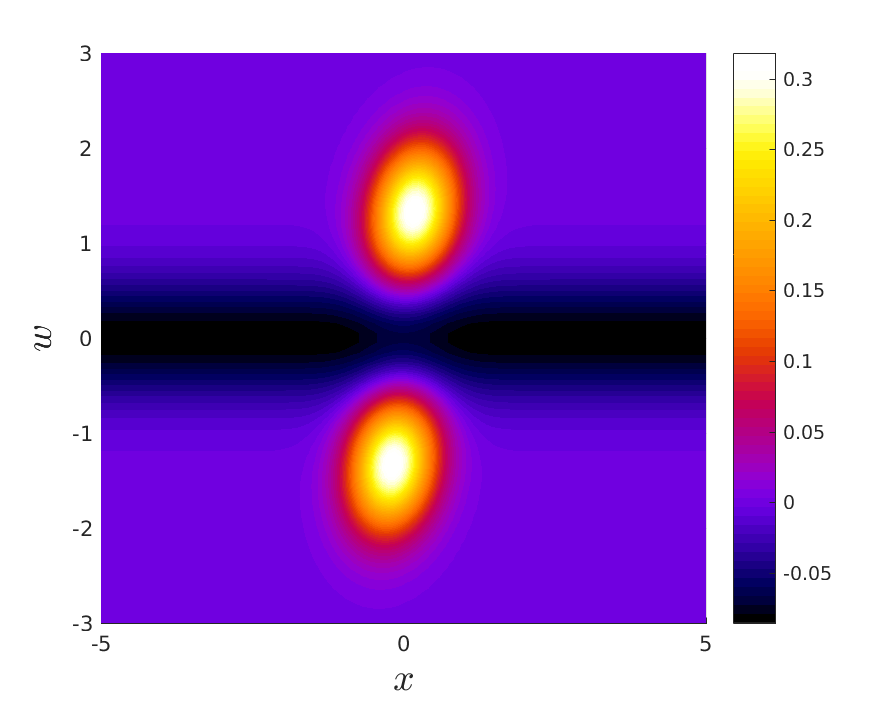}}
\subfigure[$t=3$, M$^3$C]{\includegraphics[scale=0.32]{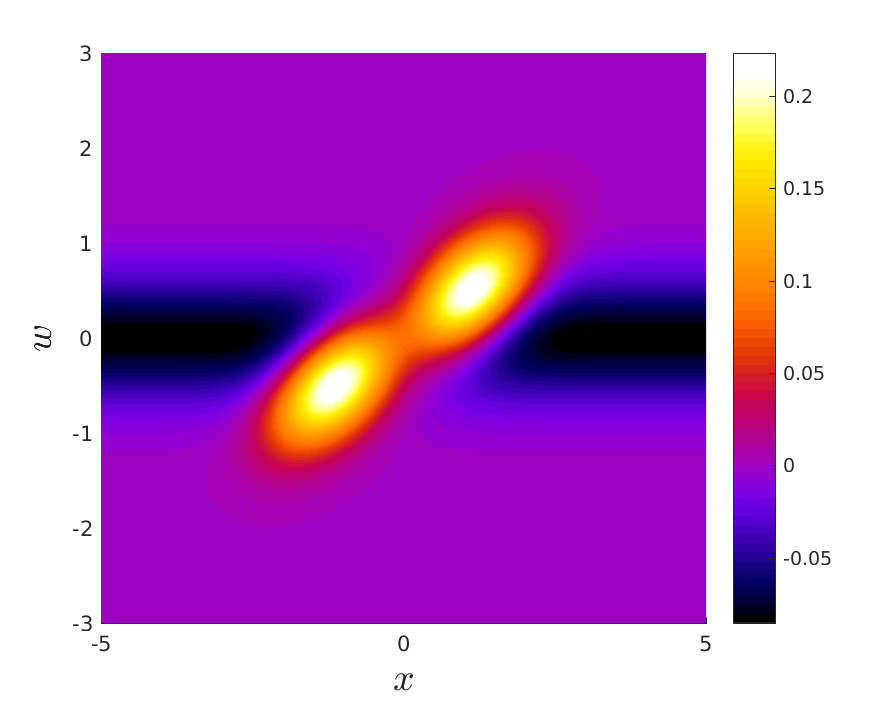}}
\subfigure[$t=6$, M$^3$C]{\includegraphics[scale=0.32]{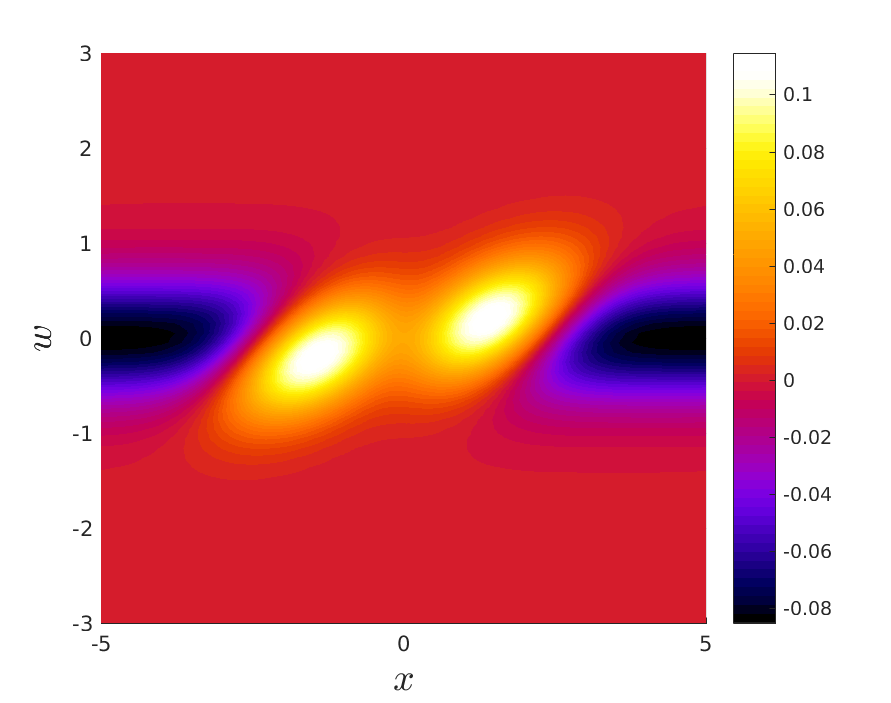}}
\caption{Example 3. Time evolution of the distribution function, expected solution over time for the MC and the M$^3$C methods. The top images report the expected solution computed with MC for $t=1$, $t=3$ and $t=6$. The bottom images
report the expected perturbation from the steady state equilibrium computed with the M$^3$C method for $t=1$, $t=3$ and $t=6$.}
\label{swarming1}
\end{figure}

\subsection{Example 3: Swarming model with uncertainties}
Finally, the last example is devoted to a Vlasov-Fokker-Planck equation describing the swarming behavior of large group of agents. It is worth to observe how for this problem one steady state solution is provided by the global Maxwellian, which is a locally stable pattern, see \cite{CFTV,DFT}. We compare the numerical solution of the problem making use of MC and M$^3$C scheme analyzed in Section \ref{sec:4}.

We consider an uncertain self--propelled swarming model described by the Vlasov-Fokker-Planck equation (\ref{eq:MF_general}) characterized by (\ref{eq:swarming_general}). This describes the time evolution of a distribution function $f(x,w,\theta,t)$ which represents the density of individuals in position $x\in\RR^{d_x}$ having velocity $w\in\RR^{d_w}$ at time $t>0$. The initial data consists in a bivariate normal distribution of the form
\begin{equation}
f_0(x,w) = C( f_0^A(x,w) + f_0^B(x,w)),
\end{equation}
where
\begin{equation}
f_0^A(x,w) = \dfrac{1}{2\pi \sqrt{\sigma_x^2 \sigma_{w}^2}}\exp\Big\{-\dfrac{1}{2}\Big( \dfrac{(x-\mu_{x})^2}{\sigma_x^2}+\dfrac{(w-\mu_{w,A})^2}{\sigma_w^2}\Big) \Big\}
\end{equation}
and
\begin{equation}
f_0^B(x,w) = \dfrac{1}{2\pi \sqrt{\sigma_x^2 \sigma_{w}^2}}\exp\Big\{-\dfrac{1}{2}\Big( \dfrac{(x-\mu_{x})^2}{\sigma_x^2}+\dfrac{(w-\mu_{w,B})^2}{\sigma_w^2}\Big) \Big\}
\end{equation}
with $\mu_x=0$, $\sigma_x=0.25$, $\mu_{w,A}=-\mu_{w,B}=1.5$, $\sigma_w^2=0.25$ and $C>0$ is a normalization constant. The uncertainty is present in the diffusion coefficient, i.e. $D=D(\theta)=0.2+0.1\theta$, and it is distributed accordingly to $\theta\sim U([-0.1,0.1])$.
\begin{figure}[t]
\centering
\subfigure[$t=1$, MC]{\includegraphics[scale=0.32]{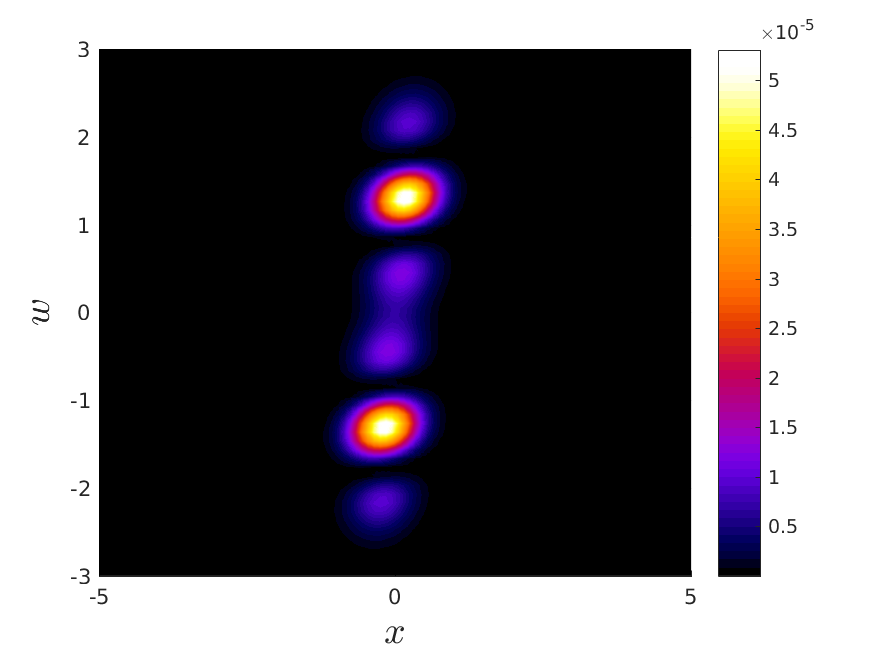}}
\subfigure[$t=3$, MC]{\includegraphics[scale=0.32]{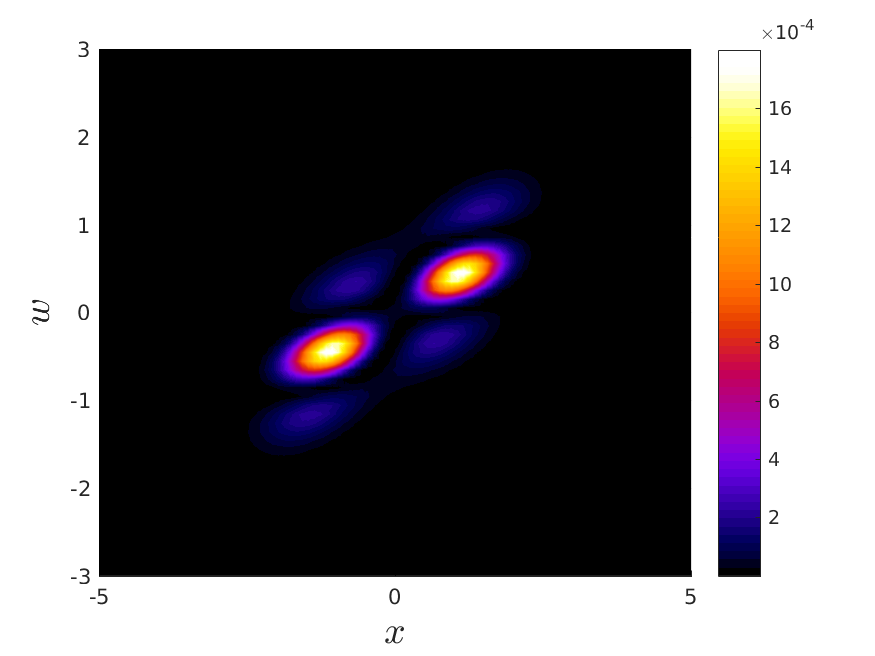}}
\subfigure[$t=6$, MC]{\includegraphics[scale=0.32]{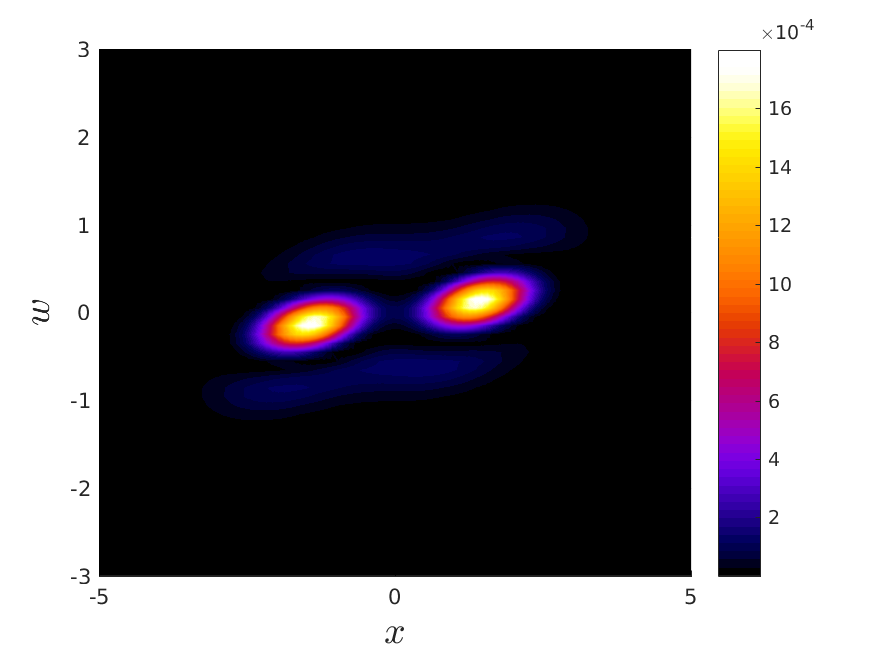}}\\
\subfigure[$t=1$, M$^3$C]{\includegraphics[scale=0.32]{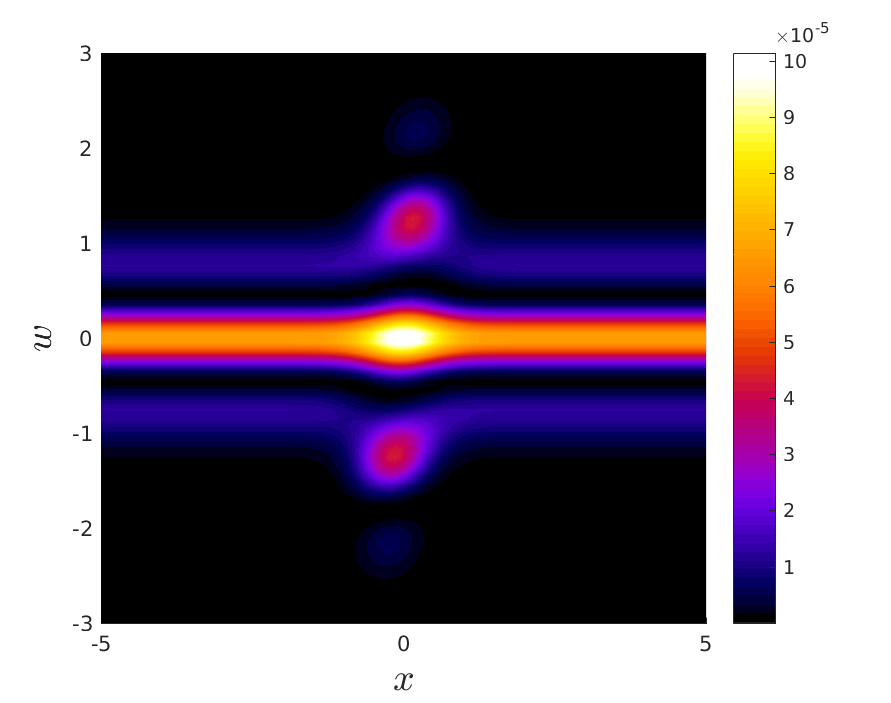}}
\subfigure[$t=3$, M$^3$C]{\includegraphics[scale=0.32]{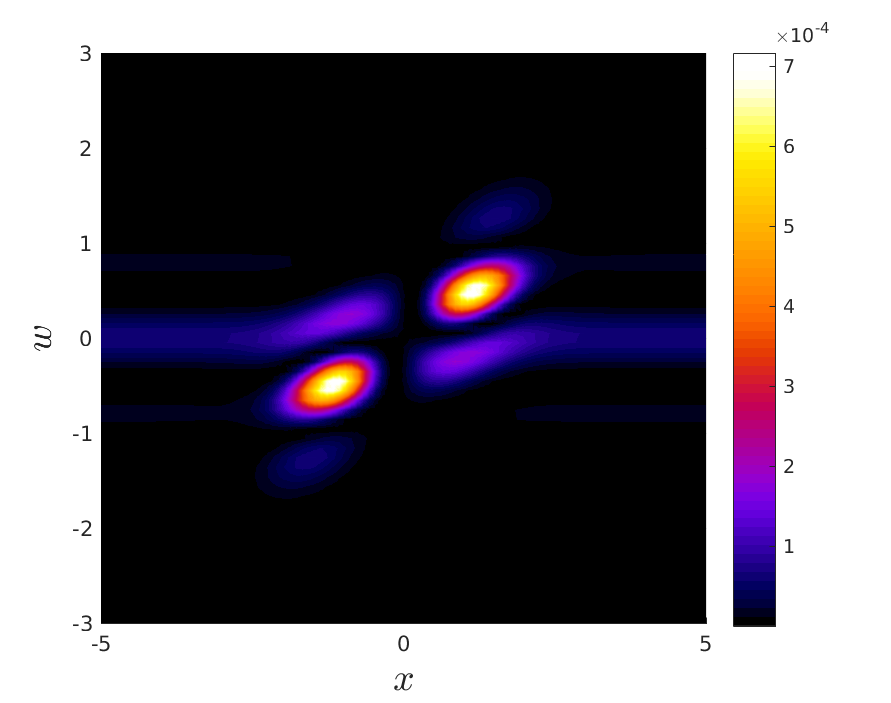}}
\subfigure[$t=6$, M$^3$C]{\includegraphics[scale=0.32]{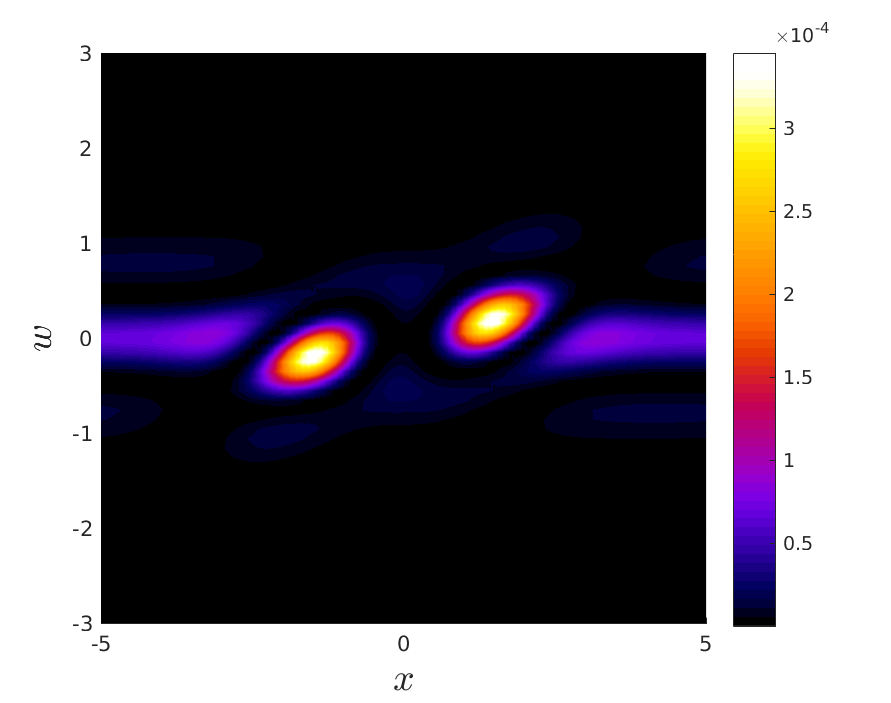}}
\caption{Example 3. Time evolution of the variance of the asymptotic solution over time for the MC and the M$^3$C methods. The top images report the variance computed with MC for $t=1$, $t=3$ and $t=6$. The bottom images report the variance of the perturbation from the steady state equilibrium computed with the M$^3$C method for $t=1$, $t=3$ and $t=6$.}
\label{swarming2}
\end{figure}

We compute the solution by using the structure preserving scheme discussed in Section \ref{sec:structure_preserving} for solving the homogeneous Fokker-Planck equation and we combine this method with a WENO scheme for the linear transport part. A second order time splitting approach joins the two discretization in space and velocity space. More in details, we compare a Monte Carlo collocation with the Micro-Macro collocation discussed in Section 3. The number of cells in space is fixed to $N_x=100$, in velocity space to $N_v=100$ while the number of random inputs is fixed to $M=50$. The solution is averaged over $10$ different realization and the final time is fixed to $T=250$. The size of the domain is $[0,L]$ with $L=10$ in space and $[-L_v,L_v]=[-3,3]$ in velocity space.
\begin{figure}[t]
\centering
\subfigure[Monte Carlo]{\includegraphics[scale=0.4]{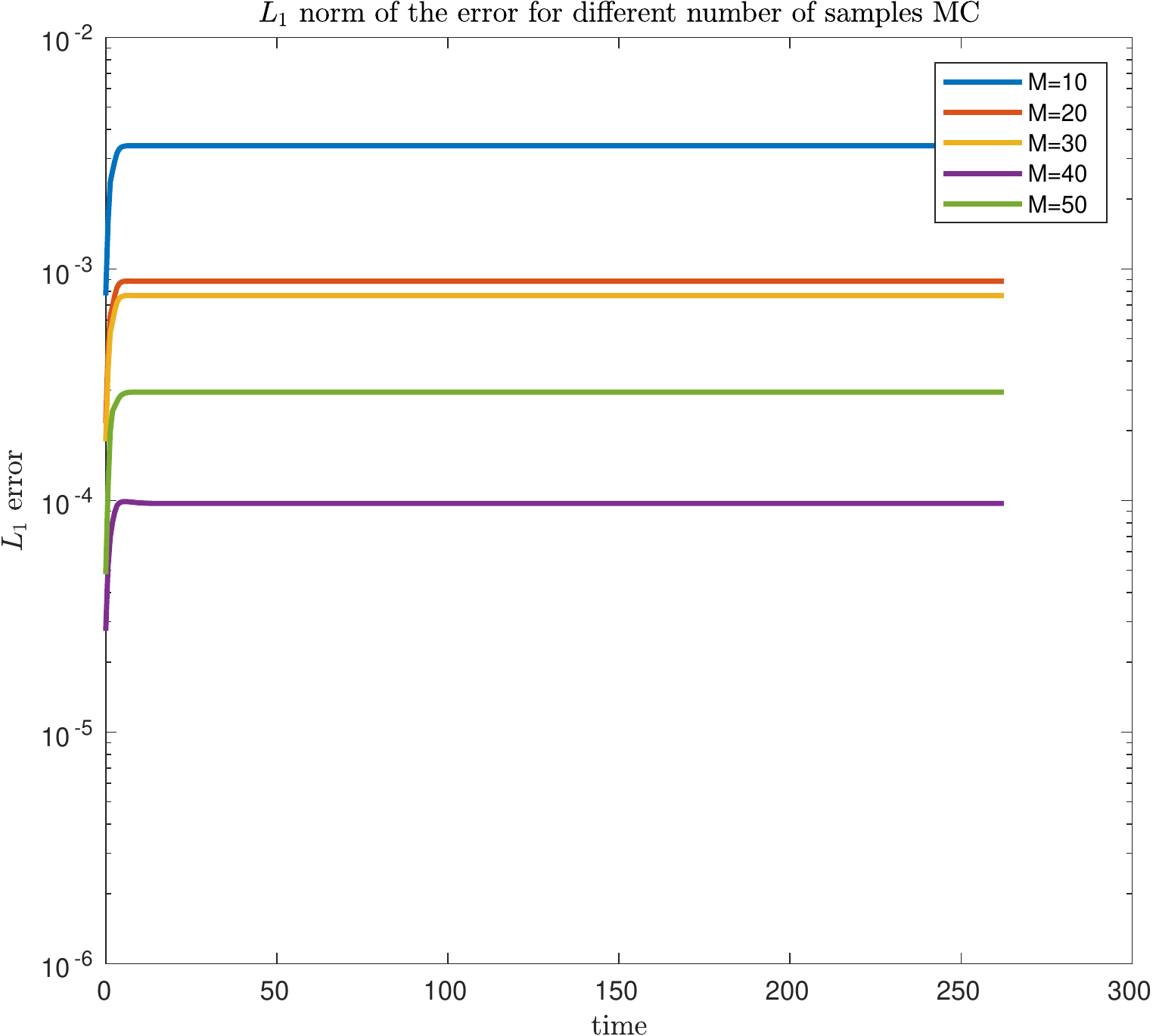}}
\subfigure[Micro-Macro Monte Carlo]{\includegraphics[scale=0.4]{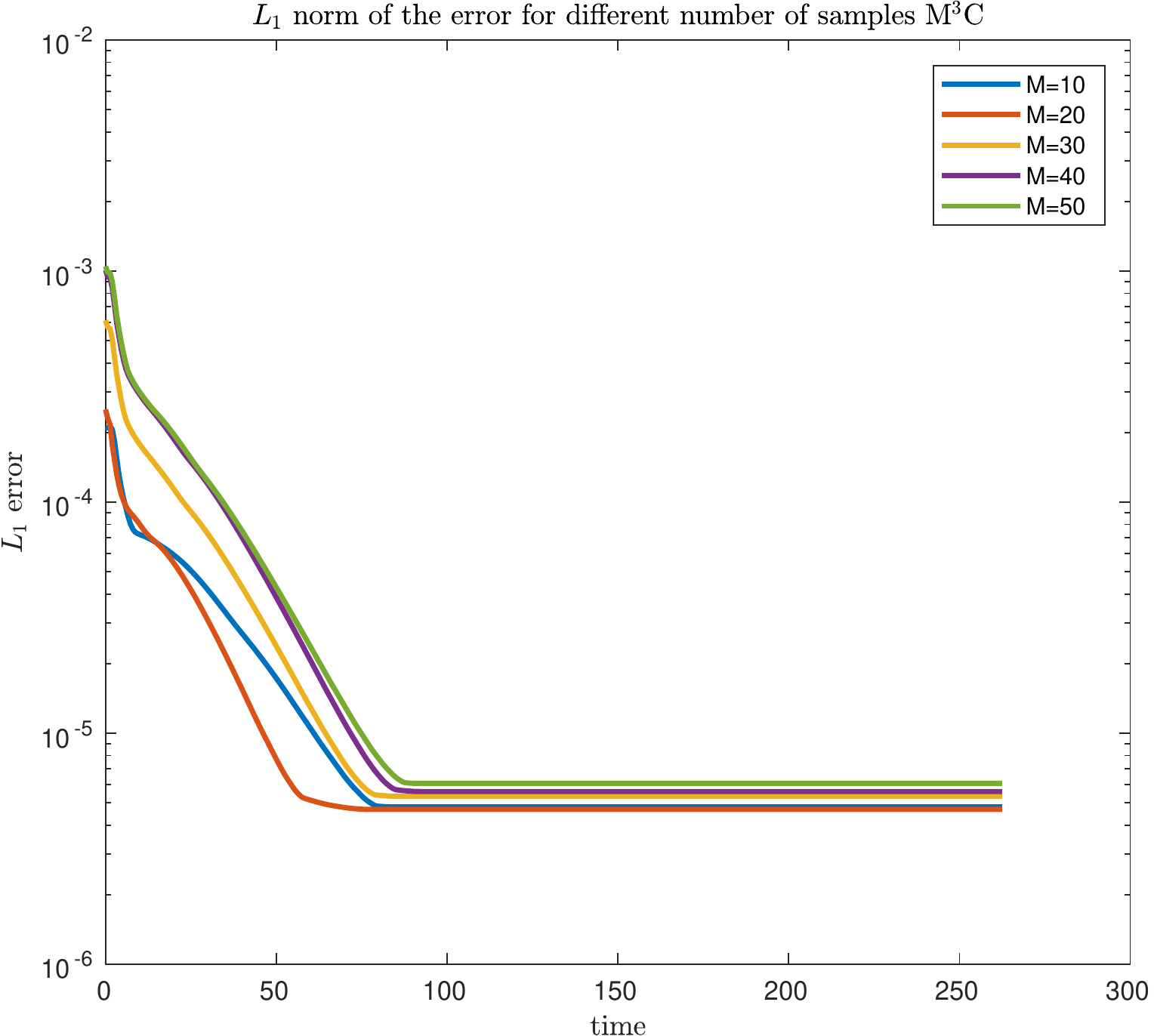}}\\
\caption{Example 3. Left: Estimation of the $L_1$ error of the expected distribution over time computed for an increasing number of 
random inputs $M$ for the MC method. Right: Estimation of the $L_1$ error of the expected distribution over time computed for an increasing number of random inputs $M$ for the M$^3$C method. The error of the $MC$ method remains constant in time while for M$^3$C method the error decreases.}
\label{swarming3}
\end{figure}
In Figure \ref{swarming1} the time evolution of the expected distribution with respect to the uncertain variable is reported for different times computed by the MC approach together  with the time evolution of the expected perturbation $g$ from the steady state solution computed with the M$^3$C method. In Figure \ref{swarming2}, the variance of the distribution over time and the variance of the perturbation $g$ over time are reported, the firsts computed by the MC method, the seconds with the M$^3$C one. Finally, in Figure \ref{swarming3}, the $L_1$ norm of the error for the MC and the M$^3$C methods are reported as a function of time for different number of random inputs. The gain in computational accuracy of the M$^3$C method is clearly evident for large times. The reference solution has been computed by a collocation method which employs the Gauss nodes as quadrature nodes with $M=100$ random inputs.

\section*{Acknowledgements}
The research that led to the present survey was partially supported by the research grant \emph{Numerical methods for uncertainty quantification in hyperbolic and kinetic equations} of the group GNCS of INdAM. MZ acknowledges support from GNCS and "Compagnia di San Paolo" (Torino, Italy).

\end{document}